\documentclass{siamltex}
\usepackage[latin1]{inputenc}% erm\"oglich die direkte Eingabe der Umlaute
\usepackage[T1]{fontenc} % das Trennen der Umlaute

\usepackage{fancyhdr}
\fancypagestyle{empty}{}%
\usepackage{graphicx,xcolor}
\usepackage{subfigure}          %% Unterabbildungen

\usepackage[english]{babel}
\usepackage[autostyle=true,german=quotes]{csquotes}
\usepackage{placeins}

\usepackage{tabularx}

\usepackage{amsmath}
\usepackage{amsfonts}
\usepackage{amssymb}
\usepackage{mathabx}

\usepackage{lmodern,textcomp}

\usepackage{algorithm}
\usepackage{algpseudocode}

\usepackage{enumitem}
\newlist{Aufz}{enumerate}{10}
\setlist[Aufz]{label*=\arabic*.}

\usepackage{xspace}
\newcommand\largeparbreak{\par\bigskip}

\newtheorem{Definition}{Definition}
\newtheorem{Lemma}{Lemma}
\newtheorem{Example}{Example}
\newcommand*\tageq{\refstepcounter{equation}\tag{\theequation}}
\newcommand{\R}{\mathbb{R}\xspace}
\newcommand{\N}{\mathbb{N}\xspace}
\newcommand{\C}{\mathbb{C}\xspace}

\newcommand{\cP}{\mathcal{P}\xspace}

\newcommand{\cU}{\mathcal{U}\xspace}

\newcommand{\cC}{\mathcal{C}\xspace}
\newcommand{\cW}{\mathcal{W}\xspace}
\newcommand{\cK}{\mathcal{K}\xspace}

\newcommand{\cO}{\mathcal{O}\xspace}

\newcommand{\PDE}{\textsf{PDE}\xspace}

\newcommand{\FE}{\textsf{FE}\xspace}

\newcommand{\BiCR}{\textsf{BiCR}\xspace}
\newcommand{\BiCGstab}{\textsf{BiCGstab}\xspace}
\newcommand{\PCR}{\textsf{PCR}\xspace}
\newcommand{\CG}{\textsf{CG}\xspace}
\newcommand{\CR}{\textsf{CR}\xspace}
\newcommand{\rhs}{\textsf{RHS}\xspace}

\newcommand{\BiCG}{\textsf{BiCG}\xspace}

\newcommand{\GMRES}{\textsf{GMRES}\xspace}

\newcommand{\GMRESkm}{\textsf{GMRES($k,m$)}\xspace}
\newcommand{\GCR}{\textsf{GCR}\xspace}
\newcommand{\GCROkm}{\textsf{GCRO($k,m$)}\xspace}
\newcommand{\MINRES}{\textsf{MINRES}\xspace}
\newcommand{\PMINRES}{\textsf{PMINRES}\xspace}
\newcommand{\RMINRES}{\textsf{R-MINRES}\xspace}
\newcommand{\SRPCRap}{\textsf{SR-PCR-ap}\xspace}

\newcommand{\KSS}{\textsf{KSS}\xspace}

\newcommand{\MVec}{\textsf{MV}\xspace}

\newcommand{\bA}{\textbf{A}\xspace}
\newcommand{\bB}{\textbf{B}\xspace}
\newcommand{\bG}{\textbf{G}\xspace}

\newcommand{\bU}{\textbf{U}\xspace}
\newcommand{\bV}{\textbf{V}\xspace}
\newcommand{\bZ}{\textbf{Z}\xspace}

\newcommand{\bP}{\textbf{P}\xspace}
\newcommand{\bPi}{\boldsymbol{\Pi}\xspace}
\newcommand{\bQ}{\textbf{Q}\xspace}
\newcommand{\bK}{\textbf{K}\xspace}

\newcommand{\tbU}{\tilde{\textbf{U}}\xspace}
\newcommand{\tbA}{\tilde{\textbf{A}}\xspace}

\renewcommand{\bf}{\textbf{f}\xspace}
\newcommand{\bx}{\textbf{x}\xspace}

\newcommand{\by}{\textbf{y}\xspace}
\newcommand{\bb}{\textbf{b}\xspace}
\newcommand{\br}{\textbf{r}\xspace}

\newcommand{\bg}{\textbf{g}\xspace}

\newcommand{\bd}{\textbf{d}\xspace}

\newcommand{\bp}{\textbf{p}\xspace}

\newcommand{\hbu}{\hat{\textbf{u}}\xspace}

\newcommand{\hbr}{\hat{\textbf{r}}\xspace}
\newcommand{\hbd}{\hat{\textbf{d}}\xspace}
\newcommand{\hbv}{\hat{\textbf{v}}\xspace}
\newcommand{\bz}{\textbf{z}\xspace}
\newcommand{\be}{\textbf{e}\xspace}
\newcommand{\bw}{\textbf{w}\xspace}
\newcommand{\bO}{\textbf{0}\xspace}

\newcommand{\bM}{\textbf{M}\xspace}

\newcommand{\bD}{\textbf{D}\xspace}
\newcommand{\bC}{\textbf{C}\xspace}

\newcommand{\bW}{\textbf{W}\xspace}

\newcommand{\bH}{\textbf{H}\xspace}
\newcommand{\bI}{\textbf{I}\xspace}
\newcommand{\bR}{\textbf{R}\xspace}
\newcommand{\bL}{\textbf{L}\xspace}
\newcommand{\bT}{\textbf{T}\xspace}

\newcommand{\bOo}{\textbf{O}\xspace}

\newcommand{\bu}{\textbf{u}\xspace}
\newcommand{\bv}{\textbf{v}\xspace}
\newcommand{\nEqns}{{n_{\text{Eqns}}}\xspace}

\newcommand{\opspan}{\textsl{span}\xspace}
\newcommand{\opImage}{\textsl{range}\xspace}

\newcommand{\obT}{\overline{\bT}\xspace}

\newcommand{\obV}{\overline{\textbf{V}}\xspace}
\newcommand{\obW}{\overline{\textbf{W}}\xspace}
\newcommand{\obD}{\overline{\textbf{D}}\xspace}

\newcommand{\obU}{\overline{\textbf{U}}\xspace}

\newcommand{\ha}{^{(1)}\xspace}
\newcommand{\hb}{^{(2)}\xspace}
\newcommand{\hc}{^{(3)}\xspace}
\newcommand{\hia}{^{(\iota)}\xspace}
\newcommand{\hib}{^{(\iota+1)}\xspace}
\newcommand{\tridiag}{\textsl{tridiag}\xspace}
\newcommand{\sign}{\textsl{sign}\xspace}

	{\end{cases}\right.}

\title{Memory-efficient Recycling of large Krylov-Subspaces for Sequences of Hermitian Linear Systems}
\author{Martin P. Neuenhofen\thanks{\texttt{Martin.Peter.Neuenhofen@rwth-aachen.de}, RWTH Aachen, Germany}
	\and
	Sven Groß\thanks{\texttt{gross@igpm.rwth-aachen.de}, Institut für Geometrie und Praktische Mathematik (IGPM),
		RWTH Aachen, Germany.}
	}

\date{\today}
\begin{document}

\maketitle

\begin{abstract}
We present a new short-recurrence residual-optimal Krylov subspace recycling \cite{PSMJM06} method for sequences of Hermitian systems of linear equations with a fixed system matrix and changing right-hand sides. Such sequences of linear systems occur while solving, e.g., discretized time-dependent partial differential equations.

With this new method it is possible to recycle large-dimensional Krylov-subspaces with smaller computational overhead and storage requirements compared to current Krylov subspace recycling methods as e.g. \RMINRES \cite{topol}.

In this paper we derive the method from the residual-optimal preconditioned conjugate residual method (\PCR, \cite[p. 182]{Saad1}) and discuss implementation issues. Numerical experiments illustrate the efficiency of our method.
\end{abstract}



%\tableofcontents

\section{Introduction}
We present a new Krylov-subspace (\KSS) method for iterative numerical solution of sequences of $\nEqns$ Hermitian linear systems
\begin{align}
	\bA  \bx\hia = \bb\hia,\qquad \iota =1,\ldots,\nEqns \,, \label{eqn:SequenceAxb}
\end{align}
with fixed $\bA \in \C^{N \times N}$, where the right-hand sides (\rhs) $\bb\hib \in \C^N$ depend on the former solution $\bx\hia$ for the \rhs $\bb\hia$. Thus the systems must be solved \emph{one after the other}. Such situations occur e.g. when applying an implicit time stepping scheme to numerically solve a non-stationary PDE. %discretization of $\partial_t u(x,t) = \partial_{xx} u(x,t)$.
Areas of application are e.g. structural dynamics \cite{PSMJM06}, topology optimization \cite{topol}, circuit analysis \cite{Circuit}, fluid dynamics \cite{NavStokes} and model reduction \cite{BeFe}.

For simplicity, in this paper we present the case where $\bA$ is regular, but from \cite{MRQLP,PCRsing} it is well-known that \CR- and \MINRES-type methods can also be applied to singular systems.

\subsection{Krylov-Subspace Methods}
Using a similar notation as in \cite{PSMJM06}, projection methods build a search space $\cU \subset \C^N$ with image $\cC = \bA  \cU$, in which an approximate solution $\bx$ for \emph{one} single system is searched, such that the \emph{residual} $\br=\bb-\bA  \bx$ is orthogonal to a test space $\cP \subset \C^N$. For \KSS methods
\begin{align}
	\cU = \bM^{-1} \cdot \cK_m(\bA  \bM^{-1};\bb) = \bM^{-1} \cdot  \operatornamewithlimits{span}_{i=0,\ldots,m-1}\lbrace (\bA  \bM^{-1})^i\,  \bb\rbrace \label{eqn:KSS_U}
\end{align}
is chosen\footnote{For an easier presentation we assume $\bx=0$ as initial guess in the introduction.}, where $\bM^{-1}$ is a cheap approximation to the inverse of $\bA$, i.e. a preconditioner (or identity when no preconditioning is used). $\cK_m$  is called \emph{Krylov-subspace} of dimension $m$. Coming from linear iterative methods or iterative refinement, this space is the most natural choice, as it adds (in case of using a reasonable preconditioner) the correction guess from the current residual to the search space for the solution, i.e. we update
\begin{align}
	\cU:=\cU + \operatorname{span}\lbrace \bM^{-1}  \br \rbrace\,. \label{eqn:KSSupdate}
\end{align}
\KSS methods augment the search space by one dimension in each iteration. For each additional dimension at least one matrix-vector-product (\MVec) with $\bM^{-1}$ and $\bA$ has to be computed, which for short recurrence methods dominate the computation time of the whole solution process. Internally \KSS methods build basis matrices $\bU,\bC,\bP \in \C^{N \times m}$ for the three spaces.

Two big goals for \KSS methods are
1. short recursions such that only a limited number of columns of the basis matrices must be stored (necessary due to hardware limitations) and also leading to lower computational costs per iteration (efficiency);
2. a good test space $\cP$ that yields small residual norms to ensure a good approximate solution w.r.t. the search space (important for convergence).
For Hermitian systems preconditioned \CR and \MINRES satisfy both conditions.

For more information about \KSS methods we refer to \cite{Saad1}.

\subsection{Krylov-Subspace Recycling}
When it comes to solving {sequences} of linear systems as in \eqref{eqn:SequenceAxb}, it seems desirable to maintain the construction idea from \eqref{eqn:KSSupdate} over the whole sequence to possibly save some of the \MVec-computations. I.e., \emph{one} high dimensional search space is successively constructed for all (or many) \rhs-es of the sequence, where always the current residual (of the current \rhs) is added to the search space. Such a construction leads to sums of \KSS-es, because \KSS information from former solution processes is reused for construction of the search space to solve the current system. Consequently such methods are called \emph{Krylov subspace recycling methods} or short \emph{recycling methods} \cite{NavStokes,PSMJM06,topol}.
\largeparbreak
In \cite{BeFe} a \GCR-based (i.e. $\bP = \bC$, $\bC$ orthogonal) \KSS recycling method is constructed by simply adding the current (preconditioned) residual as new search direction to the search space basis matrix $\bU$, cf. \cite[p. 5]{BeFe} (in their notation $\bp_k$ are the columns of $\bU$), \cite[Algo. 1]{MyReport}. This method can be characterized as a \emph{full} recycling method, as it stores the full information from all former systems. However the authors note that this method needs modifications due to memory limits and discuss remedies.

In case that the basis matrices $\bU,\bC$ have become too large to be completely stored in memory, there are two remedies known from the literature: \emph{short representations} \cite{MyReport} and basis reduction techniques (e.g. deflation \cite{Morgan1}, optimal truncation \cite{GCROT}), also summarized under the term \emph{compressings} in \cite{Circuit}.

The idea of the first approach (short representations) is to keep all the subspace information from the bases $\bU,\bC$ --- but using less storage --- by exploiting the structure of Krylov subspaces. This is what we also propose in the present paper. To our best knowledge, besides \cite{MyReport}, there is no other literature on this short representation technique.

In the other approach, based on compressings, the idea is to reduce the basis dimension by a reduction rule of the form
\begin{align}
\begin{split}
\bU_c &:= \bU\,  \boldsymbol{\eta}\\
\bC_c &:= \bC\,  \boldsymbol{\eta}\,,
\end{split} \label{eqn:Compressing}
\end{align}
where $\boldsymbol{\eta}$ has fewer columns than $\bU$ and $\bC$, leading to compressed basis matrices $\bU_c, \bC_c$. There are different approaches for the reduction, where the most common choices are to use Ritz vectors (called \emph{deflation}) \cite{PSMJM06,Morgan1} eliminating eigenvalues that slow down convergence, and \emph{optimal truncation} \cite{GCROT}.

In general in compressing methods some information is lost. So by construction the efficiency of this approach is limited compared to that of the short representation approach, as the latter is capable of conserving \emph{all} subspace information computed so far. This property only holds for the case of a constant system matrix as considered in \eqref{eqn:SequenceAxb}.

\largeparbreak

In the following we review compressing-based \KSS recycling methods from the literature.

Actually \KSS recycling methods originated from restart strategies for \GMRESkm and \GCROkm, cf. \cite{Morgan2,GMRESR,GCROT}. For both methods after $m$ iterations (of the respective \KSS method) a restart is performed. During this restart a $k$-dimensional search space basis matrix $\bU_c \in \C^{N \times k}$ ($k<m$), obtained from the iterations before by \eqref{eqn:Compressing}, is conserved and extended by a matrix $\bW\in\C^{N \times m}$ either for \GMRESkm by \emph{augmentation}
\begin{align*}
\bA  [\bU_c, \bW] &= [\bC_c,\overline{\bW}]  \begin{bmatrix}
\bI & \bB \\
\bO & \overline{\bH}
\end{bmatrix}\,,
\end{align*}
or for \GCROkm by \emph{orthogonalization}, i.e.
\begin{align}
(\bI - \bC_c  \bC_c^H) \, \bA \, \underline{ (\bI - \bC_c  \bC_c^H)^H}  \bW &= \overline{\bW} \, \overline{\bH}\,. \label{eqn:AnsatzOrtho}
\end{align}
%Details on how the matrix $\bU$ is hidden in this orthogonlization can be found in analogy in section \ref{sec:OrthKSS}.\todo{Anm. Vorkommen von $\bU$ in \eqref{eqn:AnsatzOrtho}.}
In the equations above we have $\bW = [\bw_1,\ldots,\bw_m], \obW = [\bW,\bw_{m+1}] \in \C^{N \times (m+1)}$, $\overline{\bH} \in \C^{(m+1) \times m}$ upper Hessenberg, representing the underlying Arnoldi process, and $[\bC_c,\overline{\bW}]$ \emph{orthogonal}, i.e. $[\bC_c,\obW]^H  [\bC_c,\obW] = \bI$, with $\bI$ the identity.

When not using these techniques for a restart but for recycling, the matrix $\bU_c$ in this context is called \emph{recycling basis} \cite{PSMJM06}.
\largeparbreak
In \eqref{eqn:AnsatzOrtho} the underlined factor $(\bI - \bC_c  \bC_c^H)^H$ can be dropped due to $\bC_c^H \bW=\bO$ but is written here to emphasize
%\todo{besser: emphasize statt conserve; die Tridiagonalität hat man sowieso, man braucht sie nicht künstlich zu konservieren}
the symmetry of the modified system matrix. For \RMINRES \cite{topol} this symmetry is exploited to orthogonalize $\bW$ by short recursions. As \GMRES, \MINRES, \GCR, \CR, \CG, \BiCG and \BiCR can be all derived from a common principle there are multiple publications on each of these \KSS variants for different preconditioning and compressing strategies, e.g. for \BiCG \cite{RBiCG}, \BiCGstab \cite{MaxPlanck1} and \GMRES \cite{Morgan3,HaNo}.

\subsection{Motivation and Outline}

In the following we propose a \KSS recycling method for sequences of systems \eqref{eqn:SequenceAxb}, that is capable of
\begin{enumerate}
	\item recycling a high dimensional search space $\cU$, built up by rule \eqref{eqn:KSSupdate};
	\item finding a residual-optimal\footnote{In case of preconditioning $\|\bM^{-1} \br\|_\bM$ is minimized. For details we refer to section~\ref{sec:PCR} .} solution in that space;

	\item and afterwards improving the solution by iteratively extending $\cU$ by rule \eqref{eqn:KSSupdate} only with short recurrences, and computing the residual-optimal solution in the whole space (also only using short recurrences).
\end{enumerate}
To find a residual-optimal solution in an $m$-dimensional recycled search space $\cU$ of a former system, the computational costs are as follows. Let $k,J \in \N$ with $m = k  J$ (the meaning and choice of $k,J$ will be discussed later). Our method needs $\cO(J)$ \MVec-s with $\bA$ and storage for $\cO(k)$ column vectors of length $N$ (instead of $m$ column vectors for a naive full recycling strategy). Then post-iterations are carried out where $\cU$ is extended by rule \eqref{eqn:KSSupdate} by short recurrences and the residual-optimal solution to the current \rhs is found in this extended search space, which requires 1 \MVec per post-iteration.

We call the method \SRPCRap, which stands for \emph{Short Representation based Recycling for Preconditioned Conjugate Residual with a-posteriori Optimality}.

\largeparbreak
In section~\ref{sec:SRPCR} we explain the general ideas of our method and sketch the implementation. Section~\ref{sec:num} shows numerical experiments, which illustrate the ability of our method to recycle large \KSS-es with small memory consumption. Finally, section~\ref{sec:conclusion} gives a conclusion and an outlook.

\section{\SRPCRap} \label{sec:SRPCR}
Our method is based on Preconditioned Conjugate Residual. To clarify notation and properties, we first review this method, as it is applied to solve a single system. Then we recall its properties and introduce our recycling idea.

\subsection{Preconditioned Conjugate Residual}

\subsubsection{Notation}
Throughout this text $\bU,\bD,\bV \in \C^{N \times m}$ denote basis matrices with column vectors $\bu_i,\bd_i,\bv_i \in \C^N$, $i=1,\ldots,m$. $\obU,\obD,\obV$ consist of the same basis matrices with one additional column to the right. By $\bU_{:,1:J:m}$ we denote the basis matrix that consists only of each $J^\text{th}$ column of $\bU$, starting with the first: $\bU_{:,1:J:m} = [\bu_1,\bu_{1+J},\bu_{1+2 J},\ldots] \in \C^{N \times \lfloor m/J \rfloor}$.

Let $\be_i$ denote the $i$th canonical unit vector, $\boldsymbol{1}$ a vector that contains only ones, $\bO$ a zero vector and $\bOo$ a zero matrix.
We use the symbol $\bx_0$ as an \emph{initial guess} for a solution and $\br_0 := \bb - \bA  \bx_0$ as \emph{initial residual}; for a sequence of \rhs-es we augment this notation by indices, e.g. $\bx_0\hia,\br_0\hia$ for $\iota=1, 2, \ldots$.

%For basis matrices we will use symbols like $\perp$, e.g. $\bV \perp \bD$, by which we mean properties of their ranges, i.e. $\opImage(\bV) \perp \opImage(\bD)$.

We write $\bT \in \C^{m \times m}$ for a Hermitian tridiagonal matrix. $\obT$ consists of $\bT$ with an additional row at the bottom. Its entries are denoted by
\begin{align*}
\obT &= \begin{bmatrix}
\alpha_1 	& \beta_2 	& 			& 			\\
\beta_2 	& \ddots 	& \ddots 	& 			\\
 			& \ddots 	& \ddots 	& \beta_m 	\\
 			& 			& \beta_m 	& \alpha_m 	\\
 			& 			& 			& \beta_{m+1}
\end{bmatrix} \in \C^{(m+1) \times m}\,.
%\obT &= \tridiag(\beta_2,\ldots,\beta_{m+1};\alpha_1,\ldots,\alpha_m;\overline\beta_2,\ldots,\overline\beta_m) \in \C^{(m+1) \times m}\,.
\end{align*}
%\todo{Struktur von $\obT$ und $\bT$ verständlich?}

\subsubsection{The Conjugate Residual Method (\CR)}
%For an intuition of this notation in this paper we work with
Consider an orthogonal sequence $\bv_1, \bv_2, \ldots\in\C^N$ with $\bv_i=\bA\bu_i$ obtained from a Lanczos iteration $\beta_{i+1}\bv_{i+1}= \bA\bv_i - \alpha_i\bv_i -\beta_i\bv_{i-1}$ for $i=1,\ldots,m$ and $\beta_1\equiv 0$. Equivalently, in matrix notation we obtain
\begin{align*}
\bA  \bU&= \bV \\
\bA  \bV &= \obV\, \obT
\end{align*}
with $\obV$ {orthogonal}, i.e. $\obV^H  \obV = \bI$. Starting from $\bu_1\, \parallel\, \br_0$ the \KSS $\cK_m(\bA;\br_0)$ is spanned by the columns of $\bU$. We call $\bV$ the {image} of $\bU$. \CR constructs a residual-optimal solution in $\bx_0 + \opImage(\bU)$ by subsequently orthogonalizing the residual w.r.t. each column of $\bV$ and similarly updating the numerical solution $\bx$ with the according column of $\bU$.

As $\bT$ is tridiagonal, for orthogonality each new column of $\bV$ must only be orthogonalized w.r.t. its two left neighboring columns, which makes \CR a short recurrence \KSS method. \CR is algebraically equivalent to \MINRES \cite{PaiSau75}, as both methods compute the residual-optimal solution. Orthogonality and optimality (minimal length of $\br$) are measured w.r.t. the Euclidian scalar-product $\langle \cdot , \cdot \rangle$ and norm $\|\cdot\|$.

\subsubsection{The Preconditioned Conjugate Residual Method (\PCR)} \label{sec:PCR}
Now we review preconditoning for \CR. We write $\bM \in \C^{N \times N}$ for a preconditioner with $\bM \approx \bA$, which must be Hermitian positive definite. Then instead of Euclidian measures we use $\langle \cdot , \cdot \rangle_\bM$ and $\|\cdot\|_\bM$, the $\bM$-scalar-product and its induced norm\footnote{i.e. $\langle \bu,\bv \rangle_\bM := \langle \bM  \bu,\bv \rangle$, $\|\bu\|_\bM := \sqrt{\langle \bu,\bu\rangle_\bM}$.}. The Lanczos decomposition becomes
\begin{align}
\begin{split}
\bM^{-1}  \bA  \bU &= \bM^{-1}  \bD = \bV\\
\bM^{-1}  \bA  \bV  &= \obV \, \obT
\end{split} \label{eqn:PrecLancz}
\end{align}
with $\bD=\bA\bU$, where $\obV$ is then \emph{$\bM$-orthogonal}, i.e. $\obV^H  \bM  \obV \equiv \obD^H  \obV = \bI$. $\bu_1$ is chosen parallel to $\bM^{-1}  \br_0$.

Now the \emph{preconditioned} \CR method (\PCR) does not minimize $\|\br\|$ but $\|\bM^{-1}  \br\|_\bM \equiv \sqrt{\br^H  \bM^{-1}  \br}$. For readability in the remainder, when speaking about orthogonality and optimality, these properties are always considered with respect to $\langle \cdot , \cdot \rangle_\bM$ and $\|\cdot\|_\bM$. For orthogonality w.r.t. $\langle \cdot , \cdot \rangle_\bM$ we write $\perp_\bM$. We keep on calling $\bV$ \emph{orthogonal} and the \emph{image} of $\bU$. If no preconditioner is used, i.e., $\bM = \bI$, we have $\langle \cdot , \cdot \rangle_\bM \equiv \langle \cdot , \cdot \rangle$, and $\bv_i \equiv \bd_i$, $i=1,\ldots,m+1$.

\subsubsection{Implementation, Properties and Recycling of \PCR}

\begin{algorithm}
	\caption{Preconditioned Conjugate Residual variant}
	\label{algo:PCR}
	\begin{algorithmic}[1]
	\Procedure{PCR}{$\bM,\bA,\bb,\bx_0,m$}
	\State $\tau := 1$,\quad $\bu := \bO$,\quad $\bv := \bO$
	\State $\bx := \bx_0$,\quad $\hbr := \bM^{-1}  (\bb-\bA  \bx_0)$\quad \textit{// }$\br_j = \bb-\bA  \bx_j$, $j=0,\ldots,m+1$
	\For{$j:=0,1,2,\ldots,m$} \quad \textit{// use }$\|\hbr\|_2$\textit{ or }$\theta$\textit{ as termination criterion}

		\State 	$[\hbv,\hbd,\hbu] := $\texttt{A\_func}$(\bM,\bA,\hbr)$
		\State \quad \textit{// }
					\texttt{A\_func}\textit{ is: }$\hbu := \hbr$, $\hbd:=\bA  \hbr$, $\hbv:=\bM^{-1}  \hbd$

		\State 	$\xi := \langle \hbd , \bv \rangle$

		\State 	$\bu := \hbu - \xi/\tau\,  \bu$,\quad
				$\bv := \hbv - \xi/\tau\,  \bv$

		\State  $\hat{\tau} := \langle \hbd,\bv \rangle$
		\State \quad \textit{// columns of }$\bU$\textit{: }$\bu_{j+1} := 1/\sqrt{\hat{\tau}}\,  \bu$
		\State \quad \textit{// columns of }$\bV$\textit{: }$\bv_{j+1} := 1/\sqrt{\hat{\tau}}\,  \bv$
		\State \quad \textit{// but columns of }$\bD$\textit{: }$\bd_{j+1} \neq 1/\sqrt{\hat{\tau}}\,  \hbd$\textit{\ \ - not required}

		\State \quad \textit{// estimator: }$\theta:= {\hat{\tau} + \xi^2 / \tau} \equiv \|\bM^{-1}  \bA  \br_j\|^2_\bM$

		\State 	$\alpha_j := (\tau - \xi)/\eta$,\quad
				$\beta_{j+1} := -\sqrt{\tau  \hat{\tau}}/\eta$ \quad \textit{// only if }$\obT$\textit{ is of interest}

		\State 	$\eta := \langle \hbd , \hbr \rangle$,\quad $\tau := \hat{\tau}$

		\State 	$\bx := \bx + \eta/\tau \, \bu$, \quad$\hbr := \hbr - \eta/\tau \, \bv$\quad \textit{// }$\hbr \equiv \bM^{-1}  \br_{j+1}$

	\EndFor
	\State \Return $\bx$
	\EndProcedure
	\end{algorithmic}
\end{algorithm}

Algorithm \ref{algo:PCR} gives an implementation of the \PCR method as described in section \ref{sec:PCR}.

For given $\bM,\bA,\bb,\bx_0$, \PCR constructs in $m$ iterations the search space
\begin{align}
	\cU = \bM^{-1} \cdot \cK_m(\bA  \bM^{-1};\bd) \label{eqn:SearchSpace}
\end{align}
with $\bd=\br_0$
and computes a solution $\bx \in \bx_0 + \cU$ such that $\|\bM^{-1}  \br\|_\bM$ is minimized, with $\br$ the associated residual to $\bx$. This is equivalent to
\begin{align*}
\bx &= \bx_0 + \bU  \bV^H  \br_0 \equiv \bx_0 + \bU  \bU^H  \bA  \bM^{-1}  \br_0\tageq\label{eqn:RecFormula}\\
\br &= \br_0 - \bD  \bV^H  \br_0\,,
\end{align*}
where $\bU$ is a basis of $\cU$ and $\bV$ its orthogonal image\footnote{meaning $\bV = (\bM^{-1}  \bA)  \bU$, ~$\bV^H  \bM  \bV = \bI$.}.

Eqn. \eqref{eqn:RecFormula} also provides the residual-optimal solution $\bx$ in the search space $\bx_0 + \cU$ if we replace $\cU$ in \eqref{eqn:SearchSpace} by choosing a different vector $\bd$, e.g. a residual from a former system. Consider as an example the two systems $\bA\bx\ha = \bb\ha$ and $\bA\bx\hb = \bb\hb$. Assume that we have stored the basis matrix $\bU$ of $\cU$ from the first solution process (with $\bd=\br_0\ha$). Then by virtue of formula \eqref{eqn:RecFormula} we could compute an optimal solution to the second system in $\bx_0\hb + \cU$, in which $\bx_0\hb$ is the initial guess for the second system with respective residual $\br_0\hb$.

However, in general $\bU$ does not fit into the memory. This is why in the next subsection we introduce a scheme to store $\bU$ with considerably smaller memory requirements. A discussion on how to address stability issues can be found in section~\ref{sec:OrthKSS}.

\subsection{Basic Idea: Block Krylov Matrices}
We see from Algorithm~\ref{algo:PCR} that we can get columns of $\obU$ and all the entries of $\obT$ on the fly while solving for one \rhs. The idea is now to store \emph{enough} of this data to be able to cheaply compute matrix-vector-products with $\bU$ and $\bU^H$ when using formula \eqref{eqn:RecFormula} to compute residual-optimal solutions for subsequent \rhs-es.

The key idea is the notion of a block Krylov matrix as in the following definition.
\begin{Definition}[Block Krylov Matrix]
	Let $\bZ \in \C^{N \times N}$ and $\bB \in \C^{N \times q}$ be two arbitrary matrices. For $d \in \N$ with $d \, q \leq N$ we define the \emph{block Krylov matrix} $K_d(\bZ;\bB)$ by
	\begin{align*}
		K_d(\bZ;\bB) := [\bB,\bZ  \bB,\bZ^{\,2}  \bB,\ldots,\bZ^{\,d-1}  \bB] \in \C^{N \times (d \, q)}\,.
	\end{align*}
\end{Definition}

The following lemma enables the computation of matrix-vector-products with $\bU$ and $\bU^H$ without storing the full matrix.

\begin{Lemma}[Short Representation] \label{lem:sr}
	Let $\bM,\bA \in \C^{N \times N}$, $\bU \in \C^{N \times m}$, $\bT \in \C^{m \times m}$,  as given in \eqref{eqn:PrecLancz}, with $m = k  J$, $k,J \in \N$. Define the permutation matrix $\bPi \in \C^{m \times m}$ column-wise by
	\begin{align*}
		\bPi \, \be_{1+i  J + j} = \be_{1+j  k + i},\quad j=0,\ldots,J-1,\ i=0,\ldots,k-1\,,
	\end{align*}
	and the upper triangular matrix $\bR \in \C^{m \times m}$ column-wise by
	\begin{align*}
	\bR \, \be_{1+i  J + j} = \bT^{j-1}  \be_{1+i  J},\quad j=0,\ldots,J-1,\ i=0,\ldots,k-1\,.
	\end{align*}
	Finally define $\tbU = \bU_{:,1:J:m} \in \C^{N \times k}$. Then $\bR$ has full rank and the following equality holds:
	\begin{align}
		\bU  \bR = K_J(\bM^{-1}\bA;\tbU) \, \bPi\,. \label{eqn:ShortRep}
	\end{align}
	We call $(\tbU,\bPi,\bR)$ a \emph{short representation} of $\bU$.
\end{Lemma}

This result can be found by simple calculations. $\bR$ has $\cO(m  J)$ entries and can be computed in $\cO(m  J)$. In \cite[fig. 5 \& 6]{MyReport} structures of $\bR,\bPi$ are given for sample values of $k,J$.
\largeparbreak
The proposed definition and lemma are useful in the sense, that instead of $\bU$ only the matrices $\tbU,\bR$ must be stored\footnote{The permutation represented by $\bPi$ is uniquely defined by $k,J$, cf. Lemma~\ref{lem:sr}, hence the matrix $\bPi$ does not have to be stored.}.
%\todo{Habs als Fußnote; da $\bPi$ in $\cO(m)$ ist, erscheint das nebensächlich.}
This short representation of $\bU$ enables the computation of matrix-vector-products with $\bU$ and its transpose by applying matrix-vector-products with $K_J(\bM^{-1}\bA;\tbU)$ and its transpose, respectively. The matrix-vector-products with $K_J(\bM^{-1}\bA;\tbU)$ and its transpose in turn can be computed by a Horner and a power scheme, respectively.

The idea for computation of a product of the form $\bz = K_J(\bM^{-1}\bA;\tbU) \cdot \by$, where $\by = (\by_1^T,\by_2^T,\ldots,\by_J^T)^T \in \C^m$, $\by_i \in \C^k$, is given by the Horner-type scheme
\begin{align*}
	\bz = \tbU  \by_1 + \bM^{-1}  \bA  \left( \tbU  \by_2 + \bM^{-1}  \bA  \left( \ldots \right) \right)\,,
\end{align*}
or, equivalently, algorithmically by: %\todo{Check: Danke für die Korrektur.}
\begin{algorithmic}[1]
	\State $\bz := \tbU  \tilde{\by}_{J}$
	\For{$j=J-1,\,J-2,\,\ldots,\,1$}
	\State $\bz := \bA  \bz$
	\State $\bz := \bM^{-1}  \bz$
	\State $\bz := \bz + \tbU  \tilde{\by}_{j}$
	\EndFor
\end{algorithmic}

Similarly, a product $\by = K_J(\bM^{-1}\bA;\tbU)^H \cdot \bz$, where, as above, $\by = (\by_1^T,\by_2^T,\ldots,\by_J^T)^T$, can be computed by the power scheme
\begin{align*}
\by_j = \tbU^H  \left((\bA  \bM^{-1})^{j-1} \, \bz \right)\,,\quad j=1,\ldots,J\,,
\end{align*}
or, equivalently, by:
\begin{algorithmic}[1]
	\For{$j=1,\ldots,J-1$}
	\State $\tilde{\by}_j := \tbU^H  \bz$
	\State $\bz := \bM^{-1}  \bz$
	\State $\bz := \bA  \bz$
	\EndFor
	\State $\tilde{\by}_J := \tbU^H  \bz$
\end{algorithmic}

The computational cost for applying $K_J(\bM^{-1}\bA;\tbU)$ or its transpose to a vector is $J-1$ \MVec-s with $\bA$ and $\bM^{-1}$ and $J$ matrix-vector-products with the dense $N \times k$-matrix $\tbU$.
\largeparbreak

%At that point we could do the following:
Now the recycling strategy is as follows:
We solve the first system with \PCR. During this process we collect data $\tbU$, $\bT$ to obtain a short representation of the search space $\cU$ that \PCR used. Then for a subsequent system with \rhs $\bb\hia$ we can use the short representation to recycle that search space $\cU$, i.e. by formula \eqref{eqn:RecFormula} we compute a residual-optimal solution in $\bx_0\hia + \cU$. Such a solution from a recycled search space is called \emph{recycling solution}. However, when using such a recycling strategy in practice, there are two main issues that have not been discussed yet:
\begin{itemize}
	\item \textsl{Problem A:} Due to loss of orthogonality in $\bV$ in equation \eqref{eqn:RecFormula} and stability issues with the block Krylov matrix and $\bR$ in \eqref{eqn:ShortRep} the size of $k$ and $J$ is limited in practice (depending on the system). If, for example, the size is limited to $k=10$ and $J=6$, how to recycle search spaces of more than $60$ dimensions?
	\item \textsl{Problem B:} When a recycling solution is found, how to improve it residual-optimally to satisfy a given accuracy demand?
\end{itemize}
In the next subsection we discuss one possible approach to overcome both problems.

\subsection{Building Orthogonal Krylov Subspaces} \label{sec:OrthKSS}

We first provide a simple result and then discuss its benefit.
\begin{Lemma}[\KSS-Orthogonality] \label{lem:ortho}
	From \eqref{eqn:PrecLancz} consider the last columns $\bd_{m},\bv_{m}$ of the basis matrices $\bD,\bV$. With these we define the modified system matrix
	\begin{align}
	\tbA := (\bI - \bd_{m}  \bv_{m}^H) \, \bA \, \underline{(\bI - \bv_{m}  \bd_{m}^H)}\,. \label{eqn:ModA}
	\end{align}
	Let $\br \in \C^N$ be an arbitrary vector that satisfies $\bV^H  \br = \bO$. Define its \KSS
	\begin{align}
	\cW = \bM^{\,-1} \cdot \cK_n(\tbA  \bM^{\,-1};\br)\,. \label{eqn:ModKSS}
	\end{align}
	Note that we used the modified system matrix in \eqref{eqn:ModKSS}. Then the space $\cW$ satisfies
	\begin{align}
	\bM^{\,-1}  \bA \cdot \cW \perp_\bM \bM^{\,-1}  \bA \cdot \opImage(\bU)\,. \label{eqn:OptApostCond}
	\end{align}
\end{Lemma}
For $\bM = \bI$ this result is obvious. For $\bM \neq \bI$ it follows by simple substitution of a symmetric splitting of $\bM$. As $\br \perp_\bM \bv_m$ by requirement, the underlined factor in \eqref{eqn:ModA} can be dropped, similar as in \eqref{eqn:AnsatzOrtho}.

In the following we explain the usefulness of this result: Call $\bW$ a basis of $\cW$ and $\bZ$ its orthogonal image. Further let $\bU$ and its orthogonal image $\bV$ be given by \eqref{eqn:PrecLancz}. Condition \eqref{eqn:OptApostCond} means $\bV \perp_\bM \bZ$, or $[\bV,\bZ]^H  \bM \, [\bV,\bZ] = \bI$. Thus a residual-optimal solution in $\opImage(\bU)+\opImage(\bW)$ can be computed by iteratively orthogonalizing\footnote{To be more precise: For a residual $\br$ we orthogonalize  the preconditioned residual $\hbr \equiv \bM^{-1}  \br$ onto columns of $\bV$ w.r.t. $\langle \cdot , \cdot \rangle_\bM$. Note that the computation of the $\bM$-scalar-product can be avoided using $\langle \bV,\hbr \rangle_\bM \equiv \langle \bV,\br\rangle$. In Algorithm~\ref{algo:PCR} we used $\langle \bV,\br\rangle \equiv \langle \bD,\hbr \rangle$ instead.} a residual onto columns of $[\bV,\bZ]$ and updating its numerical solution accordingly by columns of $[\bU,\bW]$.
\largeparbreak

Regarding \textsl{Problem A}, assume we used \PCR to construct a large Lanczos decomposition with basis $\bU \in \C^{N \times m}$ and tridiagonal matrix $\bT \in \C^{m \times m}$. Let us further assume that due to the large size of $\bU$ its numerical orthogonality has very low accuracy. So large errors in the recycling solution would occur if we applied \eqref{eqn:RecFormula} directly with $\bU$. Additionally let us assume that $\bT$ is ill-conditioned, so $\bR$ of $\bU$'s short representation would probably also be ill-conditioned.

To limit these influences on the recycling solution, we suggest to split the basis $\bU$ into a few \emph{basis blocks}, e.g. two. Then we have $\bU = [\bU_\text{I},\bU_\text{II}]$. The Lanczos decomposition has the property that neighboring basis vectors are orthogonal w.r.t. each other with high accuracy \cite{LancPartReortho}. So as now the blocks $\bU_i \in \C^{m_i \times m_i}$, $m_\text{I} + m_\text{II} = m$, $i \in \lbrace \text{I},\text{II}\rbrace$, are smaller than $\bU$, their respective orthogonality property\footnote{i.e. $\langle \bM^{-1}  \bA  \bU_i  \be_j , \bM^{-1}  \bA  \bU_i  \be_h \rangle_\bM = \delta_{j,h}$, $\forall i,h =1,\ldots,m_i$, $i \in \lbrace \text{I},\text{II}\rbrace$.} is more accurate. With this result, we can use \eqref{eqn:RecFormula} to orthogonalize a residual subsequently on each block (comparable to as \MINRES and \PCR do for columns).

We did not address the question yet how to keep each basis block in memory. Our solution is to store a short representation for each block. To do so we split the tridiagonal matrix
$\bT$ from above into $\bT = \diag(\bT_\text{I},\bT_\text{II}) + \textbf{E}$, where $\bT_i \in \C^{m_i \times m_i}$, $i \in \lbrace \text{I},\text{II}\rbrace$, are Hermitian tridiagonal matrices with probably better condition than $\bT$. $\textbf{E}$ consists only of $\beta_{m_\text{I}+1}$ in two subdiagonal entries. With each $\bT_i$ we can build a short representation for each $\bU_i$ of higher precision. This is done by using \eqref{eqn:ModA}, details can be found in section \ref{sec:Example} below.
\largeparbreak

Regarding \textsl{Problem B}, we consider the situation that for a subsequent system a recycling solution $\bx\hia$ has been computed with the old recycling data from \eqref{eqn:PrecLancz}, such that $\bV^H  \br\hia = \bO$. To improve the solution, we take it as initial guess for \PCR with an orthogonalization approach, where $\bA$ is replaced by the modified system matrix $\tbA$ from \eqref{eqn:ModA}, where $m$ denotes the number of columns of $\bU,\bD,\bV$ (i.e. in the formula we use the last column of each of these three matrices). As all conditions of Lemma~\ref{lem:ortho} are satisfied, \PCR will improve the solution in a residual-optimal way, as it builds a search space $\cW$ that fulfills condition \eqref{eqn:OptApostCond}.

\paragraph{Implementation}
To embed the orthogonalization approach for \PCR in algorithm \ref{algo:PCR} to address \textsl{Problem B}, \texttt{A\_func} in line 5 must be replaced by \texttt{modA\_func}:
\begin{algorithmic}[1]
\Procedure{modA\_func}{$\bM,\bA,\hbr$}
\State $[\hbv,\hbd,\hbu] := $\texttt{A\_func}($\bM,\bA,\hbr$)\quad \textit{// }$\hbr \perp \bd_m$\textit{ already}
\State $\gamma := \langle \bv_{m} , \hbd \rangle$
\State $\hbv := \hbv - \gamma \, \bv_{m}$,\quad $\hbu := \hbu - \gamma \, \bu_{m}$
\State \Return $\hbv,\hbd,\hbu$
\EndProcedure
\end{algorithmic}
We see that for the recycling the column $\bd_m$ does not need to be stored.

\subsubsection{An Explanatory Example} \label{sec:Example}

To illustrate the recycling strategy we present a fictitious example where $\nEqns = 3$ subsequent systems $\bb\ha,\bb\hb,\bb\hc$ are solved with system matrix $\bA$ and preconditioner $\bM^{-1}$.

We start by solving $\bb\ha$ with \PCR within $\hat m = 102$ iterations. During this we store the tridiagonal matrix $\bT$ up to iteration $m=90$ and from the full basis matrix $\bU$ the columns $\bu_1,\bu_6,\bu_{11},\bu_{16},\bu_{21},\ldots,\bu_{86}$ (choosing $J=5$). In addition to that we also store the columns $\bu_{30},\bv_{30}$, $\bu_{60},\bv_{60}$, and $\bu_{90},\bv_{90}$. From the stored columns we construct the three matrices $\tbU_\text{I} = [\bu_1,\bu_6,\bu_{11},\bu_{16},\bu_{21},\bu_{26}]$, $\tbU_\text{II} = [\bu_{31},\bu_{36},\bu_{41},\bu_{46},\bu_{51},\bu_{56}]$ and $\tbU_\text{III} = [\bu_{61},\bu_{66},\bu_{71},\bu_{76},\bu_{81},\bu_{86}]$. We declare the block matrices $\bU_\text{I} = [\bu_1,\bu_2,\ldots,\bu_{30}]$, $\bU_\text{II} = [\bu_{31},\bu_{32},\ldots,\bu_{60}]$ and $\bU_\text{III} = [\bu_{61},\bu_{62},\ldots,\bu_{90}]$ for later reference, but do not store them.

Choosing the parameters $k=6,J=5$, the corresponding index permutation is described by $\bPi \in \C^{30 \times 30}$ in Lemma~\ref{lem:sr}. With $\bT = \diag(\bT_\text{I},\bT_\text{II},\bT_\text{III}) + \textbf{E}$, $\bT_i \in \C^{30 \times 30}$ the diagonal blocks of $\bT \in \C^{90 \times 90}$, we construct $\bR_i \in \C^{30 \times 30}$ from input arguments $k,J,\bT_i$ for all $i \in \lbrace \text{I}, \text{II}, \text{III} \rbrace$, as described in Lemma~\ref{lem:sr}. We denote the modified system matrices\footnote{We stress that the vectors of $\bD$, that are used in these matrices, do not have to be stored, as we do not apply the modified system matrices explicitly but in a chain with the preconditioner, cf. \texttt{modA\_func}.} by
\begin{align*}
	\tbA_\text{I} &:= (\bI - \bd_{30}  \bv_{30}^H) \, \bA, \\%   (\bI - \bv_{30}  \bd_{30}^H)\\
	\tbA_\text{II} &:= (\bI - \bd_{60}  \bv_{60}^H) \, \bA, \\ %  (\bI - \bv_{60}  \bd_{60}^H)\\
	\tbA_\text{III} &:= (\bI - \bd_{90}  \bv_{90}^H) \, \bA. %  (\bI - \bv_{90}  \bd_{90}^H)\,.
\end{align*}
With these we obtain the following short representations of $\bU_\text{I},\bU_\text{II},\bU_\text{III}$:
\begin{align*}
		\bU_\text{I}  \bR_\text{I} &= K_{5}(\bM^{-1}\bA;\tbU_\text{I}) \, \bPi,\\
		\bU_\text{II} \bR_\text{II} &= K_{5}(\bM^{-1}\tbA_\text{I};\tbU_\text{II}) \, \bPi,\\
		\bU_\text{III}  \bR_\text{III} &= K_{5}(\bM^{-1}\tbA_\text{II};\tbU_\text{III}) \, \bPi\,.
\end{align*}

Next, for the solution of $\bb\hb$ we recycle the first $90$ dimensions of the former search basis $\bU$ as follows:
From the initial guess $\bx_0\hb$ we compute the associated residual $\br_0\hb$ and then apply formula \eqref{eqn:RecFormula} for $\bU_\text{I}$ to derive the improved solution $\bx\hb_\text{I}$. For $\bx_\text{I}\hb$ in turn we compute the associated residual $\br_\text{I}\hb$ and then apply formula \eqref{eqn:RecFormula} for $\bU_\text{II}$ to derive the improved solution $\bx\hb_\text{II}$. Finally for $\bx_\text{II}\hb$ we compute the associated residual $\br_\text{II}\hb$ and then apply formula \eqref{eqn:RecFormula} for $\bU_\text{III}$ to obtain the recycling solution $\bx\hb_\text{III}$ of the recycling basis $[\bU_\text{I},\bU_\text{II},\bU_\text{III}]$. For each block $2  J$ \MVec-s are necessary for application of formula \eqref{eqn:RecFormula}. Thus we needed $30$ \MVec-s so far and storage for $24$ column vectors of length $N$, but found a residual-optimal solution in a $90$-dimensional space.

Now we want to perform a-posteriori iterations to improve $\bx\hb_\text{III}$. For that purpose we run \PCR on the system $\bA  \bx\hb = \bb\hb$ with initial guess $\bx\hb_\text{III}$ and preconditioner $\bM^{-1}$, where \texttt{A\_func} in \PCR is replaced by \texttt{modA\_func} for $\tbA_\text{III}$. After $m = 42$ iterations \PCR has found a residual-optimal solution in the whole search space of $90 + 42 = 132$ dimensions with sufficiently small residual. The basis of this $42$-dimensional search space is denoted by $\bU_\text{IV} = [\bu_{91},\bu_{92},\ldots,\bu_{132}]$.

In total, to compute the solution for $\bb\hb$ residual-optimally in a $132$-dimensional search space, only $24$ column-vectors needed to be stored (except resources during \PCR) and $30+42 = 72$ \MVec-s with $\bA$ and $\bM^{-1}$ needed to be computed.

The solution to system $\bb\hc$ can be computed in the same way as for system $\bb\hb$, by simply reusing the recycling basis blocks $[\bU_\text{I},\bU_\text{II},\bU_\text{III}]$ by short representations, respectively, by subsequently orthogonalizing a residual of a numerical solution w.r.t. each block. As a remark: Alternatively a short representation of $\bU_\text{IV}$ could be fetched during the a-posteriori iterations for system $\bb\hb$ and added to the other three recycled blocks. Then by virtue of \eqref{eqn:RecFormula} for system $\bb\hc$ a recycling solution $\bx\hc_\text{IV}$ (with a residual that is orthogonal to the images of all four recycling basis blocks) can be computed. However for the a-posteriori-iterations to improve $\bx\hc_\text{IV}$ in an optimal sense we need then the modified system matrix
\begin{align*}
	\tbA_\text{IV} = (\bI - [\bd_{90},\bd_{132}]  [\bv_{90},\bv_{132}]^H) \, \bA \, \underline{(\bI - [\bv_{90},\bv_{132}]  [\bd_{90},\bd_{132}]^H)}\,,
\end{align*}
because $\bu_{90}$ did not originate from the \KSS where $\bu_{132}$ originated from. (The underlined factor can be dropped.)

\subsection{Computational Cost}

For simplicity we assume that we recycle a \KSS of only one system (i.e. we neglect the remark above), by $\ell$ basis blocks of identical width $m = k  J$, each of which is shortly represented by a dense $N \times k$-matrix. (For the above example we had $\ell = 3$, $k=6$, $J=5$.) We neglect storage and computational cost in $\cO(\ell  m  J)$ (with $J \leq m$).

To orthogonalize a residual w.r.t. one block, one needs to compute $2  J$ \MVec-s with both $\bA$ and $\bM^{-1}$, $2  J$ scalar-products, $6  J$ AXPYs  and $J$ matrix-vector-products with both the dense $N \times k$-matrix and its transpose.

In total, for recycling a $\ell  k  J$-dimensional search space, $\ell  (k+2)$ column-vectors need to be stored. The $+2$ occurs due to the two vectors $\bu_{m_i},\bv_{m_i}$ that must be additionally stored for each block $i=1,\ldots,\ell$ for the modified system matrix.

To construct a residual-optimal recycling solution in the recycled $\ell  k  J$-dimensional search space, the cost is dominated by $2  \ell  J$ \MVec-s with $\bA$ and $\bM^{-1}$, and $\ell  J$ matrix-vector-products with a dense $N \times k$-matrix and its transpose, respectively.

\section{Numerical Experiments} \label{sec:num}

In this section we present four test cases in which we compare the proposed method \SRPCRap to the well-known method \MINRES \cite{PaiSau75}, as it is implemented in Matlab. All test cases are from partial differential equations (\PDE) discretizations. The first two are taken from the Florida Sparse Matrix Collection \cite{FSMC}, the second two stem from finite element discretizations.

\paragraph{The \enquote{Right} Right-Hand Sides}
In \cite[sec. 3.7]{Saad1} the author states that the choice of \rhs-es is not crucial for a test case. We disagree in that point as for matrices representing discrete elliptic differential operators (which are of practical relevance) it is well-known that for high-frequent solutions a \KSS solver needs more iterations to resolve these frequencies.

It also seems to be a common practice to define the \rhs by the (non-smooth) image of a smooth solution vector, e.g. in \cite{IDRstab}. However, when solving in \cite[sec. 6.4 (c)]{IDRstab} for the \rhs $u$ (smooth) instead of $F$, the according solution vector has high frequencies (which is physically meaningful, cf. to solutions of Turing bifurcation problems). To solve for this \rhs $u$, \GMRES and \MINRES would need twice as many iterations as for $F$. To sum it up, choosing a \rhs from a smooth solution vector may provide convergence curves that differ much from those for a \rhs coming from practical applications.

In addition to that, symmetry both of the domain and the \rhs-function can influence the convergence speed extremely. As an example, for $-\Delta u = 1$, $u\vert_{\partial \Omega}=0$ on a Cartesian isotropic grid of the domain $\Omega = (0,1)^2$, one can expect at least $\sqrt{8}$ times faster convergence than for a random \rhs due to symmetry (which of course also affects round-off errors and thus stability, and due to lower dimension of the full \KSS also the superlinear convergence).

To avoid any confusion about initial guesses and whether it is useful to choose the initial guess of the next system as the solution to the former system, the initial guess $\bx_0 = \bO$ is chosen for each \rhs and ten \rhs-es $\bb\ha,\ldots,\bb^{(10)}$ are constructed for each test case, that are orthogonal w.r.t. each other.

For our test cases there is additionally the question how to choose a series of \rhs-es: When e.g. choosing two \rhs-es $\bb\ha$, $\bb\hb$, where $\bb\hb$ is orthogonal to the image of the former search space\footnote{i.e. $\bb\hb \perp_\bM \bA  \bM^{-1} \cdot \cK_m(\bA  \bM^{-1};\bb\ha)$} then any \KSS recycling method is by construction just as good as its respective method without recycling, cf. example~\ref{expl:OrthRecSpace}. So we have to ensure that the recycled search space of $\bb\ha$ contains good solution candidates for the subsequent \rhs-es (\emph{utility}) and on the other hand we have to ask which kind of sequences occur in practice (\emph{practical relevance}). We make three suggestions, going from practical relevance to utility. For each sequence the \rhs-es are chosen orthogonal w.r.t. each other and independent of $\nEqns$, the number of \rhs-es, using an orthogonal\footnote{w.r.t. the Euclidian scalar-product} basis of a \KSS to some \emph{starting vector} $\bd$:
\begin{enumerate}
	\item \textsl{Sequence A:} When the next \rhs depends on a linear combination of former solutions, then the choice
	\begin{align}
		\opspan\lbrace\bb\ha,\ldots,\bb^{(q)} \rbrace = \cK_q(\bA^{-1};\bd)\quad \forall\,q=1,\ldots,\nEqns \label{eqn:RHSnatural}
	\end{align}
	seems reasonable. For Hermitian positive definite systems we know by the convergence results from \cite[eqn. 6.107]{Saad1},\cite[chap. 5.2]{GR11}, that $\cK_m(\bA;\bb\ha)$ contains good solution candidates for the subsequent \rhs-es of this sequence for sufficiently large $m$. However, this is not clear for the search space with preconditioning $\bM^{-1} \cdot \cK_m(\bA  \bM^{-1};\bb\ha)$.

	\item \textsl{Sequence B:} We generalize \eqref{eqn:RHSnatural} under preconditioning.
	\begin{align}
	\opspan\lbrace\bb\ha,\ldots,\bb^{(q)} \rbrace = \cK_q(\bM  \bA^{-1};\bM^{-1}  \bd)\quad \forall\,q=1,\ldots,\nEqns \label{eqn:RHSqnatural}
	\end{align}
	results in \eqref{eqn:RHSnatural}, if $\bM$ is splitted (e.g. Cholesky) and substituted into $\bA$. For the same reasons as above it is known that $\bM^{-1} \cdot \cK_m(\bA  \bM^{-1};\bb\ha)$ contains good solution candidates for the subsequent \rhs-es of this sequence for sufficiently large $m$.
	%Such sequences often occur in finite element (\FE) discretizations of non-stationary \PDE-s, where $\bM$ is a mass matrix.
	\item \textsl{Sequence C:} This last sequence is defined for systems where strong preconditioning is required to be able to find an iterative solution at all. Numerical experiments indicate that the nearer the preconditioner is to $\bA^{-1}$, the more suited the recycling space is \emph{only} to $\bb\ha$, cf. results for sequence A in section~\ref{sec:numPOISSON}.
	Therefor we consider the following test sequence, where the influence of the preconditioner on the utility is attenuated:
	\begin{align}
	\opspan\lbrace\bb\ha,\ldots,\bb^{(q)} \rbrace = \cK_q(\bA  \bM^{-1};\bd)\quad \forall\,q=1,\ldots,\nEqns\,. \label{eqn:RHSunnatural}
	\end{align}
	The search space $\bM^{-1} \cdot \cK_m(\bA  \bM^{-1};\bb\ha)$ often provides good solution candidates for these \rhs-es, but it does not contain the corresponding exact solution. In this setting recycling could be probably replaced by parallel solves for the different \rhs-es, as there is no dependence of a \rhs on a former solution (as no $\bA^{-1}$ occurs).
\end{enumerate}
\largeparbreak
\begin{Example}[Limits of Recycling] \label{expl:OrthRecSpace}
	We stress that there exist sequences of linear systems for which recycling is completely useless. As an example, let $N \in 2 \N$, $\bA = 1/(N+1)^2 \cdot \tridiag(-1,2,-1) \in \R^{N \times N}$, $\bM = I$, and consider a sequence of $2$ \rhs-es $\bb\ha = (\textbf{1}^T,\textbf{1}^T)^T$, $\bb\hb = (-\textbf{1}^T,\textbf{1}^T)^T$ with $\textbf{1} \in \R^{N/2}$.
	%As initial guess we take $\bx_0 = \bO$ for both systems.
	The equations can be interpreted as discretizations of Poisson problems $-\Delta u = f\hia$ in $\Omega = (-1,1)$, $u(-1)=u(1)=0$ with \rhs-functions $f\ha(x) = 1$, $f\hb(x) = \sign(x)$.

	Obviously for the first system the \rhs, the solution as well as the vectors $\bA^j  \bb\ha$ $\forall j \in \N$ are symmetric w.r.t. the symmetry axis $x=0$. In contrast, for the second system the \rhs, the solution and the vectors $\bA^j  \bb\hb$ $\forall j \in \N$ are antisymmetric w.r.t. $x=0$. Consequently, the recycling spaces $\cU\hia := \bM^{-1} \cdot \cK_N(\bA  \bM^{-1};\bb\hia)$, $\iota=1,2$, have the following properties:
	\begin{align*}
		\cU\hb  &\perp_\bM \cU\ha\,,\\
		\bM^{-1}  \bA \cdot \cU\hb &\perp_\bM \bM^{-1}  \bA \cdot \cU\ha\,,
	\end{align*}
	which means that 1.) the recycling space $\cU\ha$ is perpendicular to $\bx\hb$ and 2.) its image is orthogonal to $\bb\hb$.
	Thus for the second system a recycling solution in $\cU\ha$ does not improve the initial guess $\bx_0 = \bO$, regardless of whether a residual-optimal or an error-optimal recycling method would be used.
\end{Example}

\FloatBarrier

\subsection{CURLCURL0}

In this test case \cite{CURL} only the indefinite symmetric matrix $\bA \in \R^{11083 \times 11083}$ is given. For sake of simplicity we chose $\bd = \bA \, \textbf{1}$.

\begin{figure}
	\centering
	\includegraphics[width=0.85\linewidth]{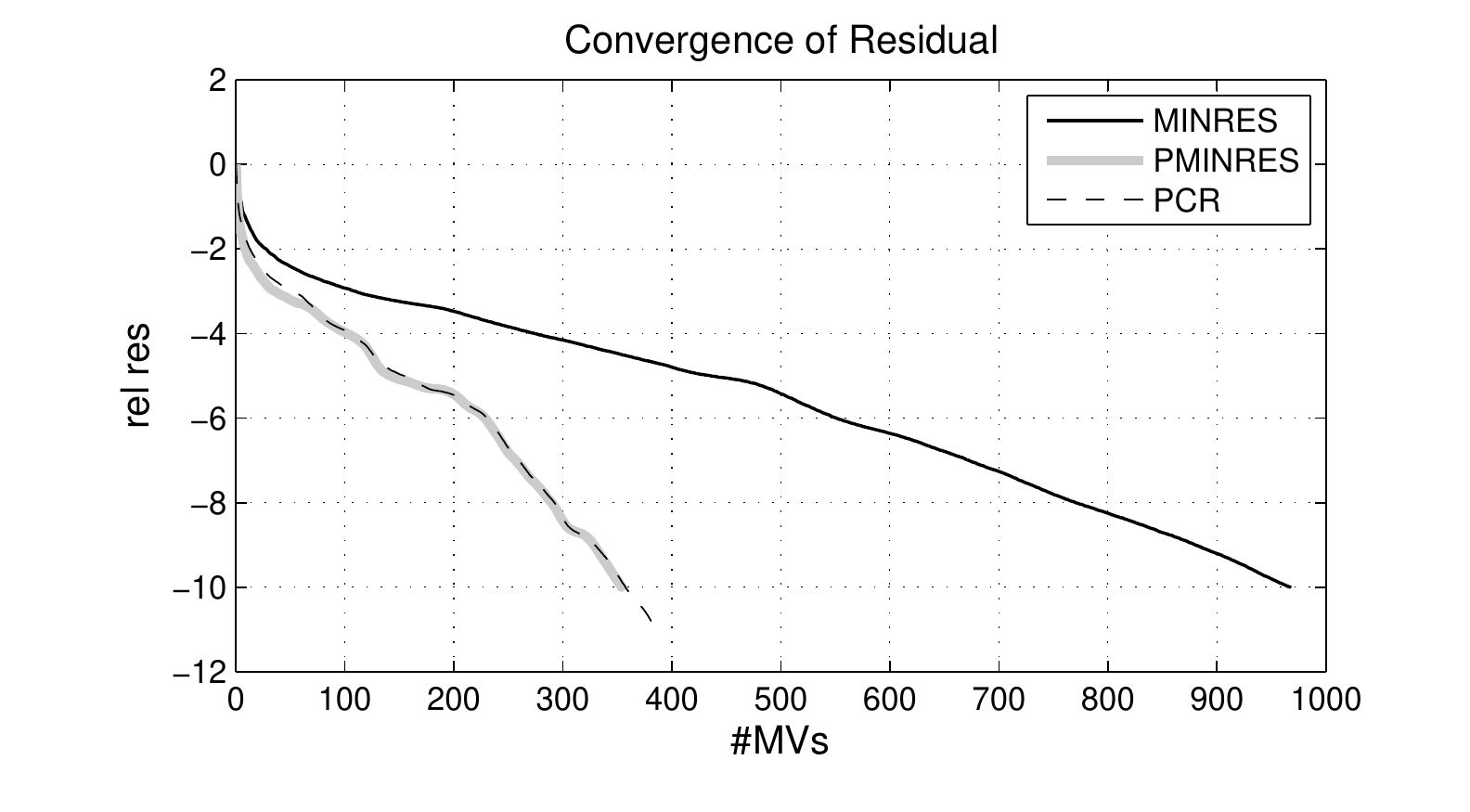}
	\caption{Convergence of \MINRES,\PMINRES and \PCR for CURLCURL0 with \rhs $\bb = \bA \, \boldsymbol{1}$.}
	\label{fig:Conv_CURL}
\end{figure}

\paragraph{System without preconditioning}
We first study the unpreconditioned system (i.e., $\bM = \bI$) with system matrix $\bA$ and the \rhs $\bb = \bd$. The convergence of $\|\br\|_2/\|\bb\|_2$ for \MINRES is shown in black in figure \ref{fig:Conv_CURL}.

\begin{figure}
	\centering
	\includegraphics[width=0.9\linewidth]{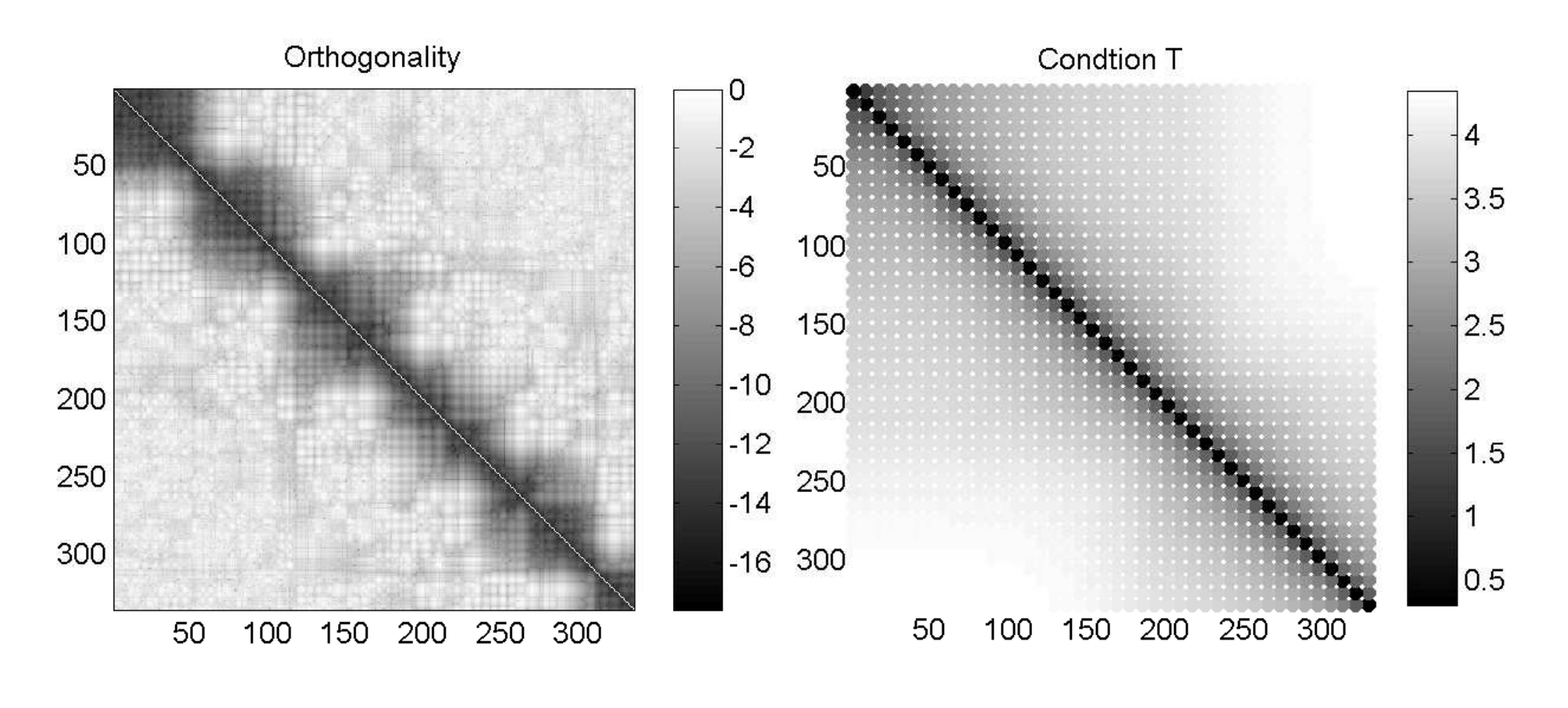}
	\caption{Matrices (left) $\bQ$ and (right) $\bG$ for unpreconditioned CURLCURL0.}
	\label{fig:Comp_CURL}
\end{figure}

Next, we want to estimate an appropriate block size for recycling, cf. section~\ref{sec:OrthKSS}. Therefore the matrix $\bQ := \bV^H  \bM \, \bV$ is computed with $\bV$ from \eqref{eqn:PrecLancz} of the solution process with \PCR. $\bQ$ should be the identity but due to round-off errors it also has non-zero off-diagonal entries. This affects the stability of the Gram-Schmidt orthogonalization that is hidden in \eqref{eqn:RecFormula}. The $\log_{10}|\,\cdot\,|$ for each entry of $\bQ$ is shown in the left of figure \ref{fig:Comp_CURL}. Right in the figure $\log_{10}|\,\cdot\,|$ is given for each entry of a symmetric matrix $\bG$ that has the entries $g_{i,j} = \kappa_2( \bT_{i:j,i:j} )$, i.e. the $2$-condition of a section of the tridiagonal matrix from \eqref{eqn:PrecLancz}. We discussed in section \ref{sec:OrthKSS} that the ill-conditioning of these sub-matrices affects the accuracy of the short representations. So as a rule of thumb, the bandwidth in which the entries of both matrices of figure \ref{fig:Comp_CURL} are small limits the block size for each basis block. These two plots give an estimate how large each basis block can be chosen. Taking larger blocks is beneficial as increasing the number of blocks leads to an increased memory consumption.

From $\bQ$ in fig.~\ref{fig:Comp_CURL} we see that e.g. a block of the first 50 basis columns could be recycled (as these are nearly orthogonal w.r.t. each other), then maybe a small block from columns 51 to 60 could be recycled, and then a block from columns 60 to 100, and so on. However this seems complicated and not so efficient, as we would need many small blocks, which needs more storage and computational effort than a few large blocks. We now turn to a study of the block size for the preconditioned system.

\paragraph{Preconditioned System}
We will see that by use of preconditioning the orthogonality and conditioning properties are improved, i.e., the off-diagonal elements of $\bQ$ and $\bG$ become smaller. Here the preconditioner
\begin{align}
\bM = \sign(\diag(\bA)) \cdot \tridiag(\bA) \label{eqn:PrecM}
\end{align}
is used which is Hermitian positive definite. The convergence of $\|\br\|_2/\|\bb\|_2$ for \PMINRES (i.e. preconditioned \MINRES with $\bM$) is given in fig.~\ref{fig:Conv_CURL} in thick gray. The iterates are algebraically equivalent to those of \PCR, but as \PCR does not iteratively compute $\br$, we plot $\|\bM^{-1}  \br\|_2/\|\bM^{-1}  \bb\|_2$ instead (dashed curve in fig.~\ref{fig:Conv_CURL}). In the remainder we assume both residual measures to be comparable.

\begin{figure}
	\centering
	\includegraphics[width=0.9\linewidth]{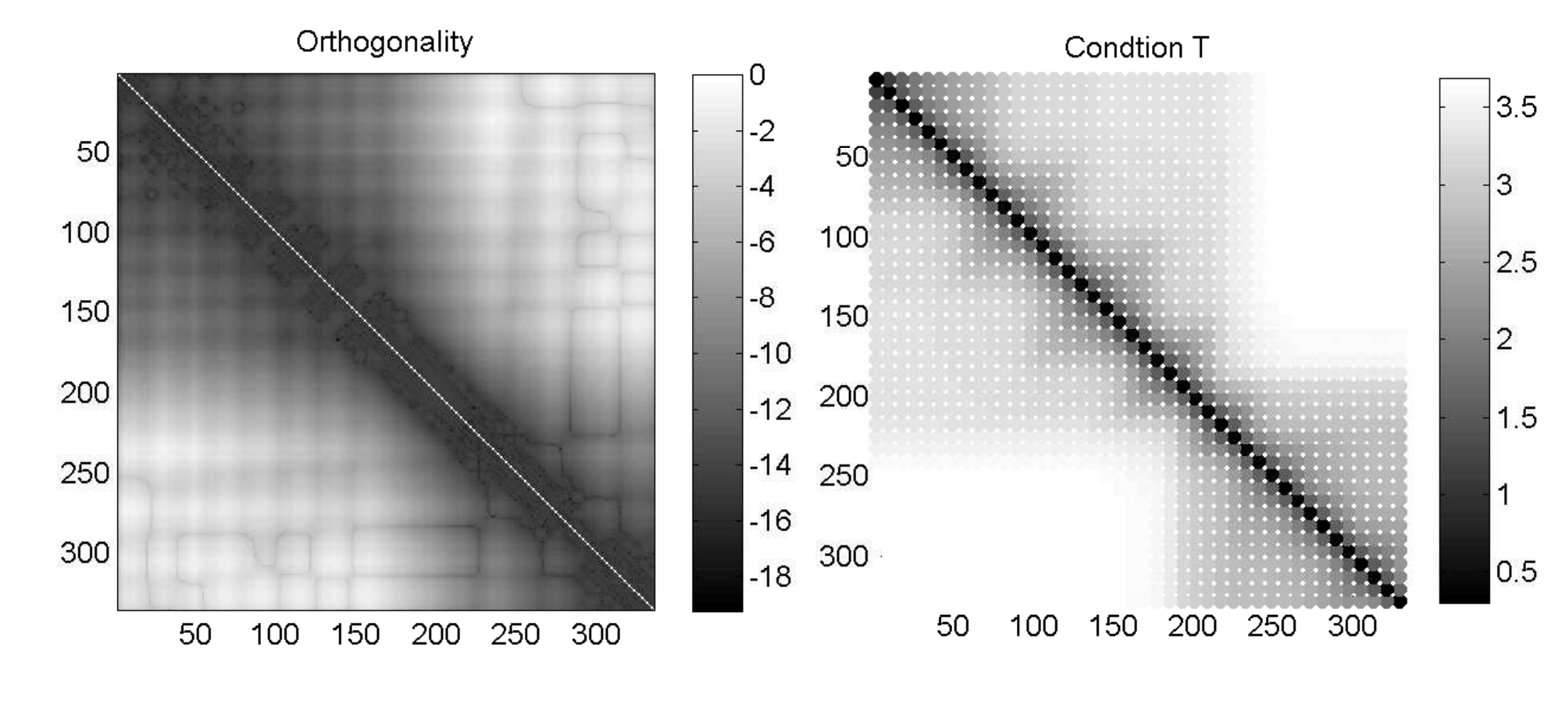}
	\caption{Matrices (left) $\bQ$ and (right) $\bG$ for preconditioned CURLCURL0.}
	\label{fig:Comp_CURLprec}
\end{figure}

Extracting the basis columns from the solution process of \PCR, the matrices $\bQ,\bG$ are shown in fig.~\ref{fig:Comp_CURLprec}. Comparing the result to fig.~\ref{fig:Comp_CURL}, we see that now the basis blocks can be chosen larger, e.g. the first 150 columns are nearly orthogonal w.r.t. each other. However, the conditioning of $\bT$'s sub-blocks still limits the block sizes. As for practical problems the matrices $\bQ,\bG$ cannot be computed, we will make very defensive choices for the block sizes in the following.

\paragraph{Solving Sequence C}
For $\bM,\bA,\bd$ we construct $\bb\ha,\bb\hb,\ldots,\bb^{(10)}$ from {sequence C}.
To solve the sequence, we first solve for $\bb\ha$ with \PCR and fetch short representations for its basis blocks. For this we choose $\ell = 7$ subsequent basis blocks of each $k=8$ vectors and $J=6$ (yielding a block size of $k  J = 48$ for each block). By this the search space $\bM^{-1} \cdot \cK_{336}(\bA  \bM^{-1};\bb\ha)$ can be recycled with a cost of $70$ stored columns and $84$ \MVec-s.
Afterwards we solve for $\bb\ha,\bb\hb,\ldots,\bb^{(10)}$ with \SRPCRap, i.e. we use formula \eqref{eqn:RecFormula} subsequently for each block, followed by residual-optimal post-iterations.

\begin{figure}
	\subfigure[]{%
		\includegraphics[width=0.48\linewidth]{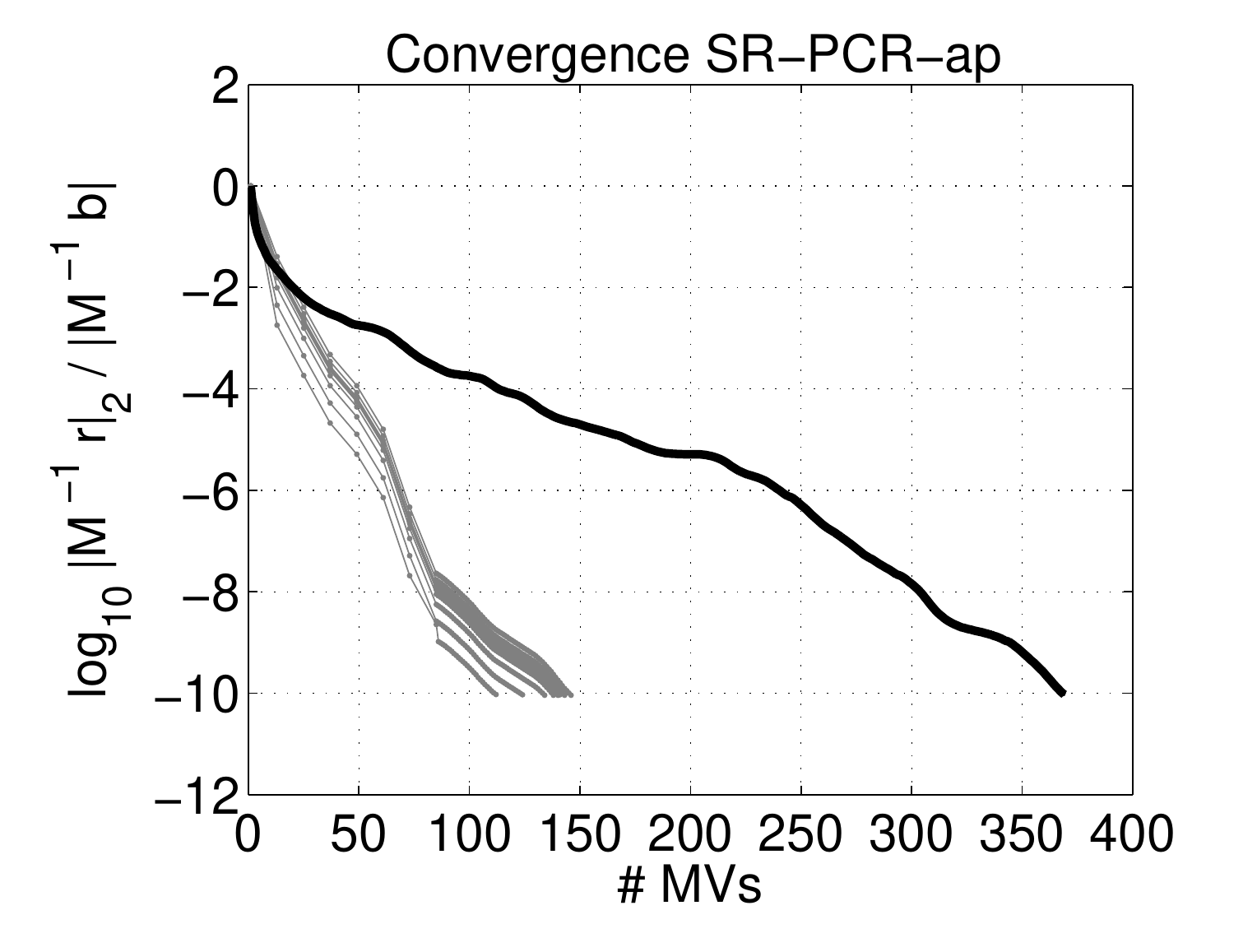}
	}
	\hfill
	\subfigure[]{%
		\includegraphics[width=0.48\linewidth]{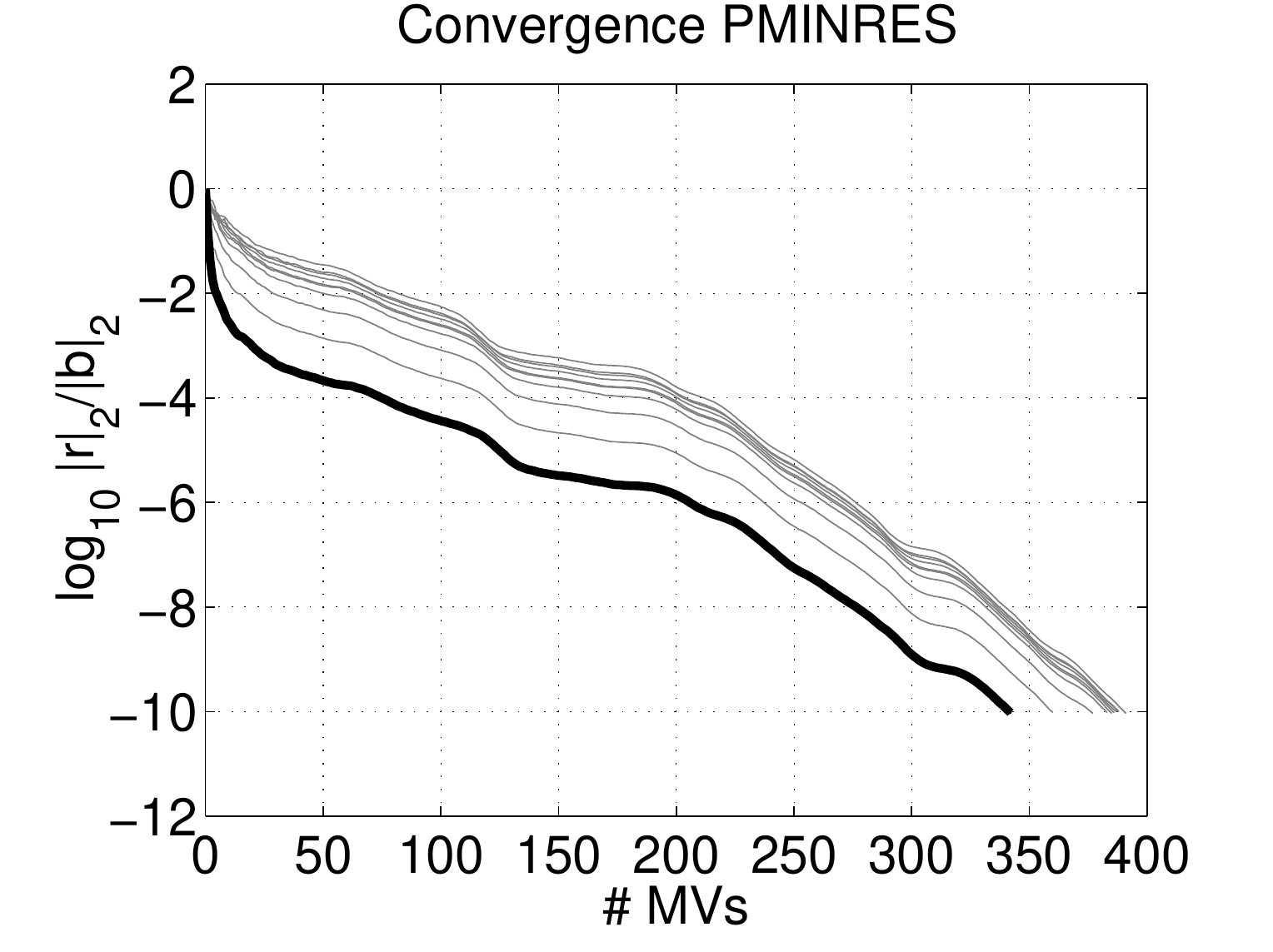}
	}
	\caption{Convergence for all \rhs-es of {sequence C} of CURLCURL0 for (a) \SRPCRap, (b) \PMINRES; thick black = first \rhs, thin gray = subsequent \rhs-es.}
	\label{fig:Prec_Conv_CURL}
\end{figure}

Fig. \ref{fig:Prec_Conv_CURL} shows the convergence of \SRPCRap (left) in comparison to \PMINRES (right). One can observe that recycling leads to a dramatic acceleration of convergence speed in this example. E.g., for a required accuracy of $10^{-8}$ \SRPCRap needs only $1/3$ of the \MVec-s compared to \PMINRES.

\subsection{SHERMAN1}
This test set \cite{Boeing} provides a symmetric indefinite matrix $\bA \in \R^{1000 \times 1000}$, $\kappa_2(\bA) \approx 2.3\cdot 10^4$, with one vector that we use for $\bd$, the starting vector of the \rhs-sequence.

We directly go for the preconditioned system using a preconditioner $\bM$ as in \eqref{eqn:PrecM} for CURLCURL0. For the \rhs-es we used {sequence B}.

As for CURLCURL0, we solve for $\bb\ha$ with \PCR, collect short representations, and then recycle these basis information using the proposed method \SRPCRap. As parameters we choose $\ell = 2$ blocks of each $k=8$ stored columns and $J=7$, and thus recycle the search space $\bM^{-1} \cdot \cK_{112}(\bA  \bM^{-1};\bb\ha)$ storing $20$ columns and applying 28 \MVec-s.

\begin{figure}
	\subfigure[]{%
		\includegraphics[width=0.48\linewidth]{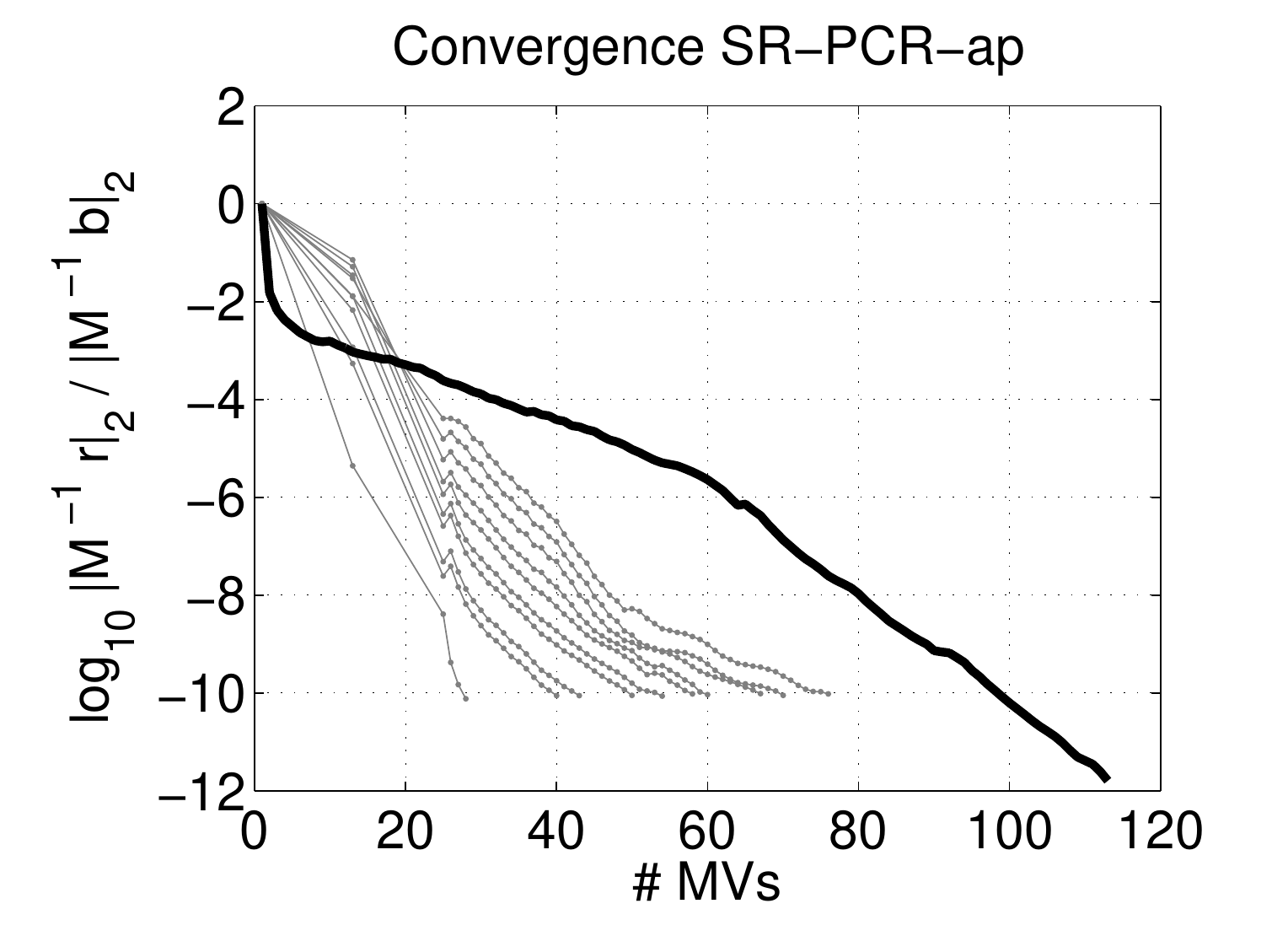}
	}
	\hfill
	\subfigure[]{%
		\includegraphics[width=0.48\linewidth]{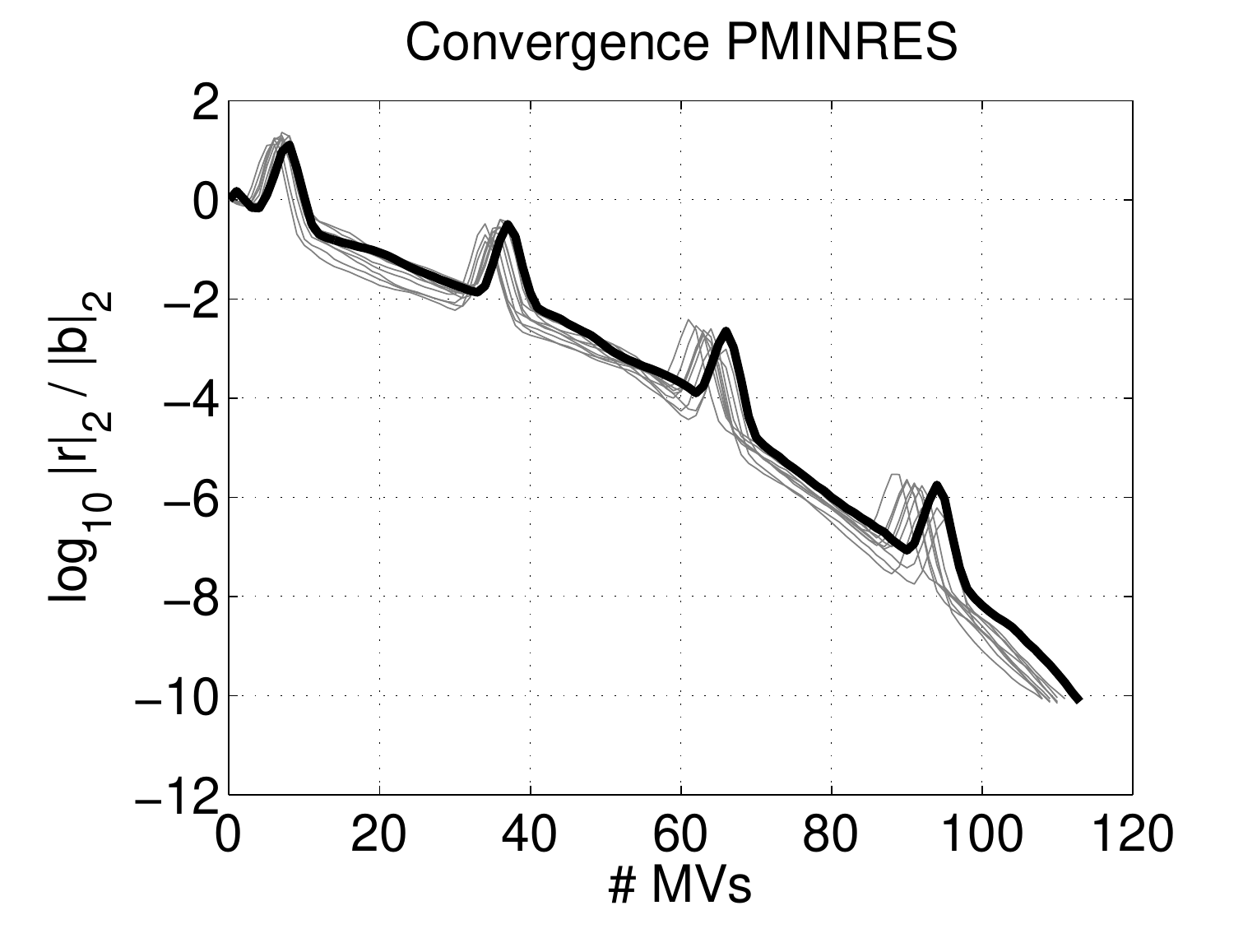}
	}
	\caption{Convergence for all \rhs-es of {sequence B} of SHERMAN1 for (a) \SRPCRap, (b) \PMINRES; thick black = first \rhs, thin gray = subsequent \rhs-es.}
	\label{fig:Prec_Conv_SH1}
\end{figure}

Fig.~\ref{fig:Prec_Conv_SH1} shows a comparison of the convergence of the residuals for \SRPCRap with those of \PMINRES. We see that by use of the recycling information the proposed method converges faster for each \rhs, but for later \rhs-es the number of \MVec-s grows, as the recycling space becomes out of date.

\subsection{Finite Element Test Cases}
\begin{figure}
	\subfigure[]{%
		\includegraphics[width=0.48\linewidth]{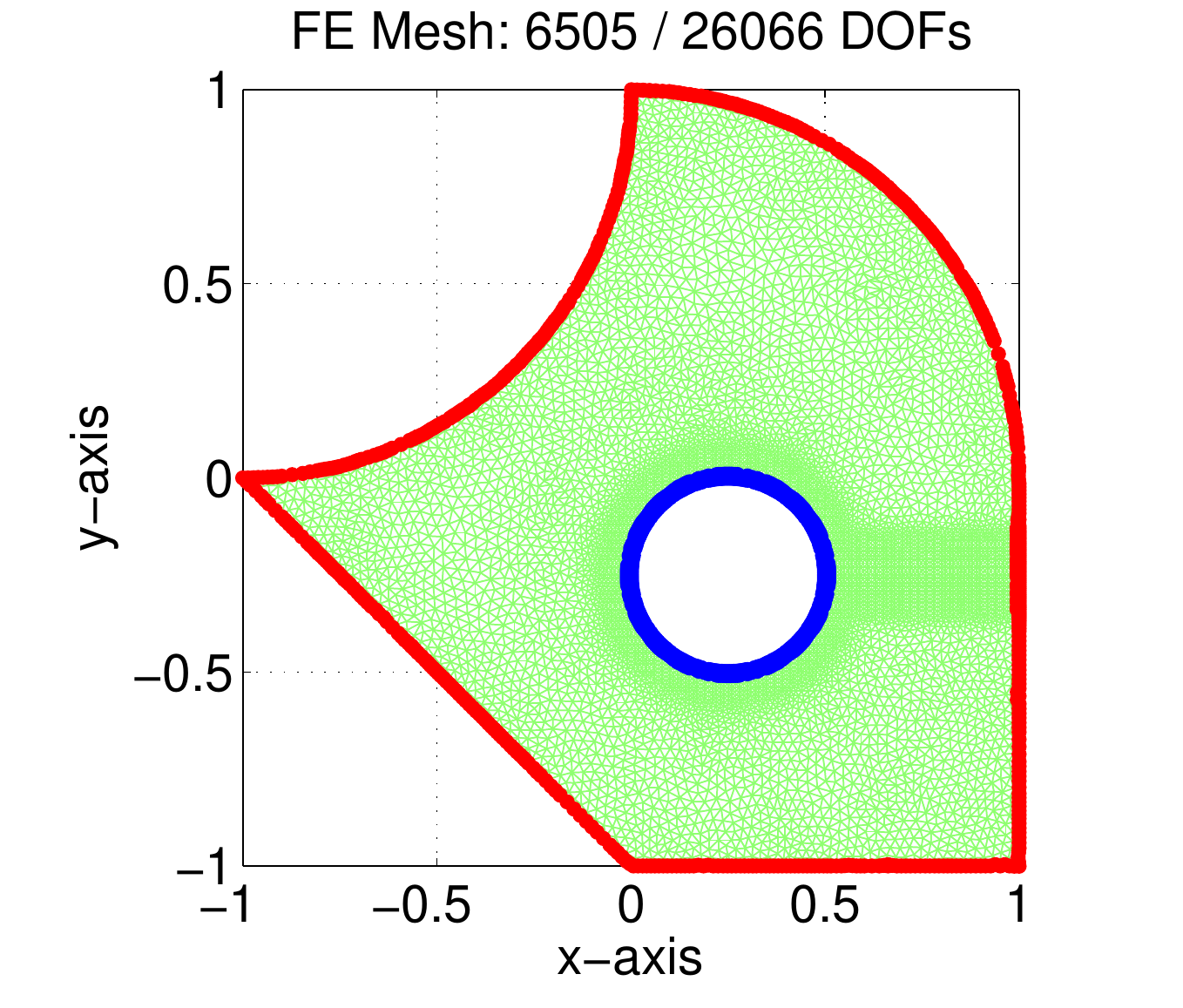}
	}
	\hfill
	\subfigure[]{%
		\includegraphics[width=0.48\linewidth]{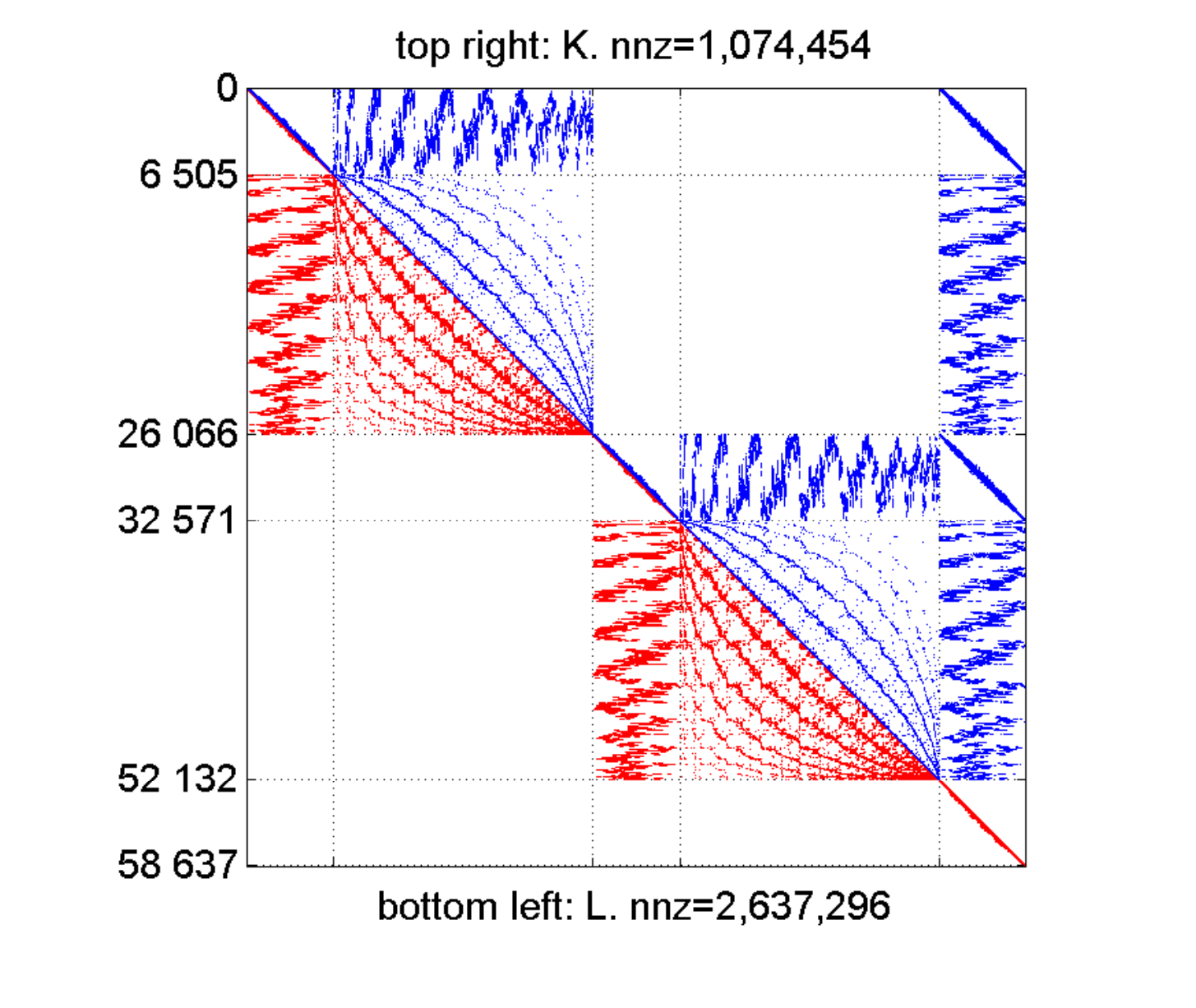}
	}
	\caption{(a) Finite Element Mesh with nodes: red = Dirichlet boundary nodes, blue = Neumann boundary nodes (inner hole). (b) Pattern of Stokes Matrix (blue, only upper right part) and its Preconditioner (red).}
	\label{fig:Mesh1}
\end{figure}

The two following test cases stem from \FE discretizations of two different \PDE-s in 2D, namely a Poisson and Stokes problem. The \FE discretizations are carried out using SimpleFEM \cite{SimpleFEM}.
The two-dimensional domain $\Omega$ and the corresponding triangular mesh are shown in figure~\ref{fig:Mesh1} (a).
%To define our test systems in a simple but realisic manner we solve both on one two-dimensional domain $\Omega$ on the tetrahedral \FE discretization of piecewise linear / quadratic shape functions given in figure \ref{fig:Mesh1} (a).
The Neumann boundary part $\Gamma_N$ is given in thick gray (inner boundary; hole), the Dirichlet part $\Gamma_D$ in black. On this mesh, piecewise quadratic \FE are chosen for the discretization of the Poisson problem. For the Stokes problem piecewise linear and piecewise quadratic \FE are used for the discretization of the pressure and the velocity field, respectively.

\subsubsection{Poisson Problem} \label{sec:numPOISSON}
We compute a numerical solution $\phi$ to
\begin{align*}
			-\Delta \phi &= 0  \quad \text{in } \Omega\\
			\vec{\nabla}\phi \boldsymbol{\cdot} \vec{n} &= 0  \quad \text{on } \Gamma_N \\
			\phi &= g_D\quad \text{on } \Gamma_D
			\qquad \text{with }
			g_D(x,y):= \begin{cases} 1, & \mbox{if } x = 1 \\ 0, & \mbox{else}  \end{cases}
\end{align*}
using piecewise quadratic \FE on the mesh shown in fig.~\ref{fig:Mesh1} (a). The resulting stiffness matrix $\bA \in \R^{26066 \times 26066}$ is symmetric positive definite with $\kappa_1(\bA) \approx 6.5 \cdot 10^4$. The non-zero \rhs $\bd$ contains the contributions from the non-homogeneous Dirichlet  boundary data. We choose the preconditioner $\bM = \bL  \bL^H$, where $\bL$ is the incomplete Cholesky decomposition \cite[Algo. 10.7]{Saad1} of $\bA$.

We try $\ell = 6$ blocks, each of dimension $48$, by choosing $k=8$, $J=6$ for each block. By this we recycle the search space $\bM^{-1} \cdot \cK_{288}(\bA  \bM^{-1};\bb\ha)$ by storing $60$ columns and performing $72$ \MVec-s.

\begin{figure}
	\subfigure[]{%
		\includegraphics[width=0.48\linewidth]{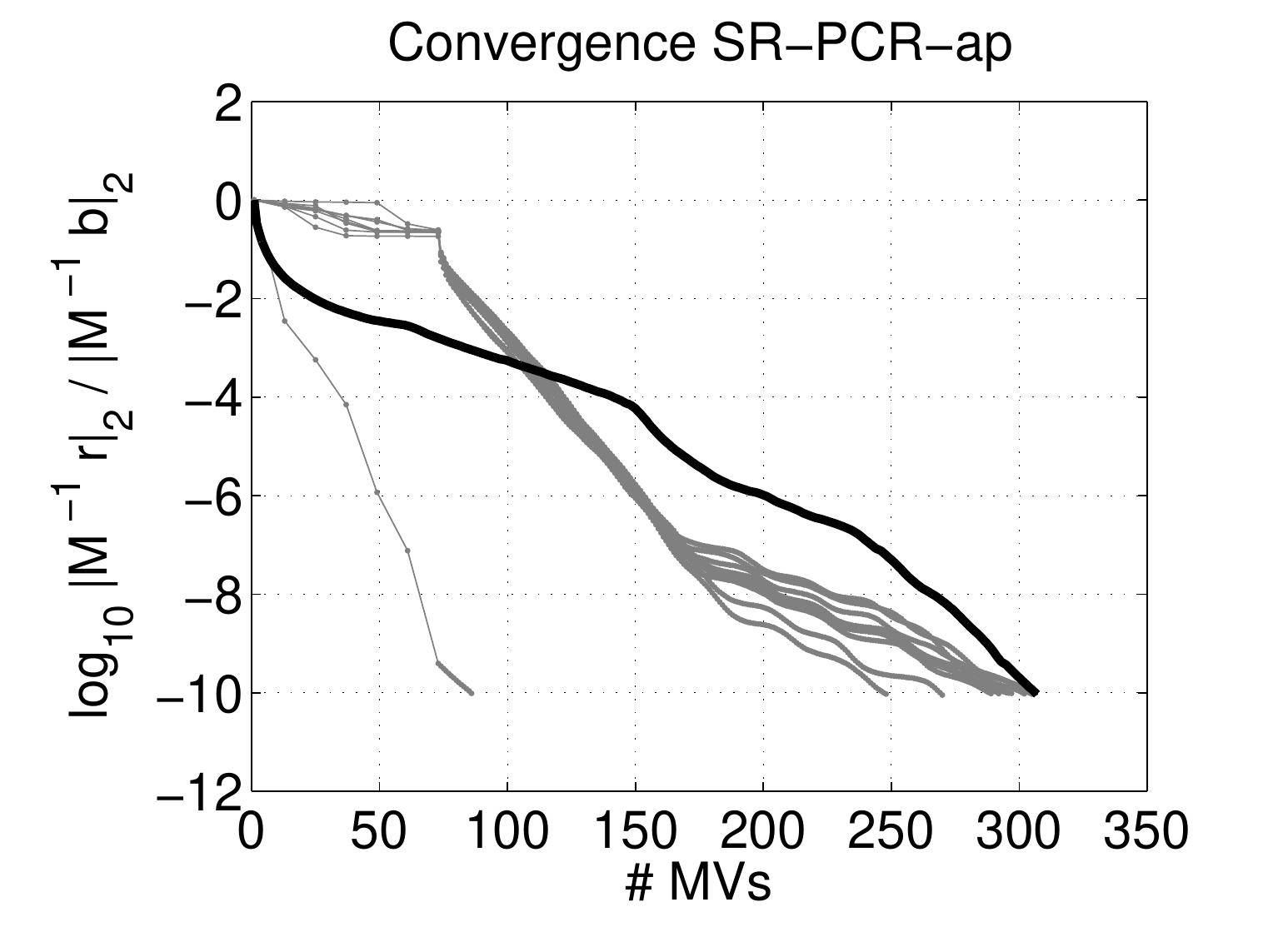}
	}
	\hfill
	\subfigure[]{%
		\includegraphics[width=0.48\linewidth]{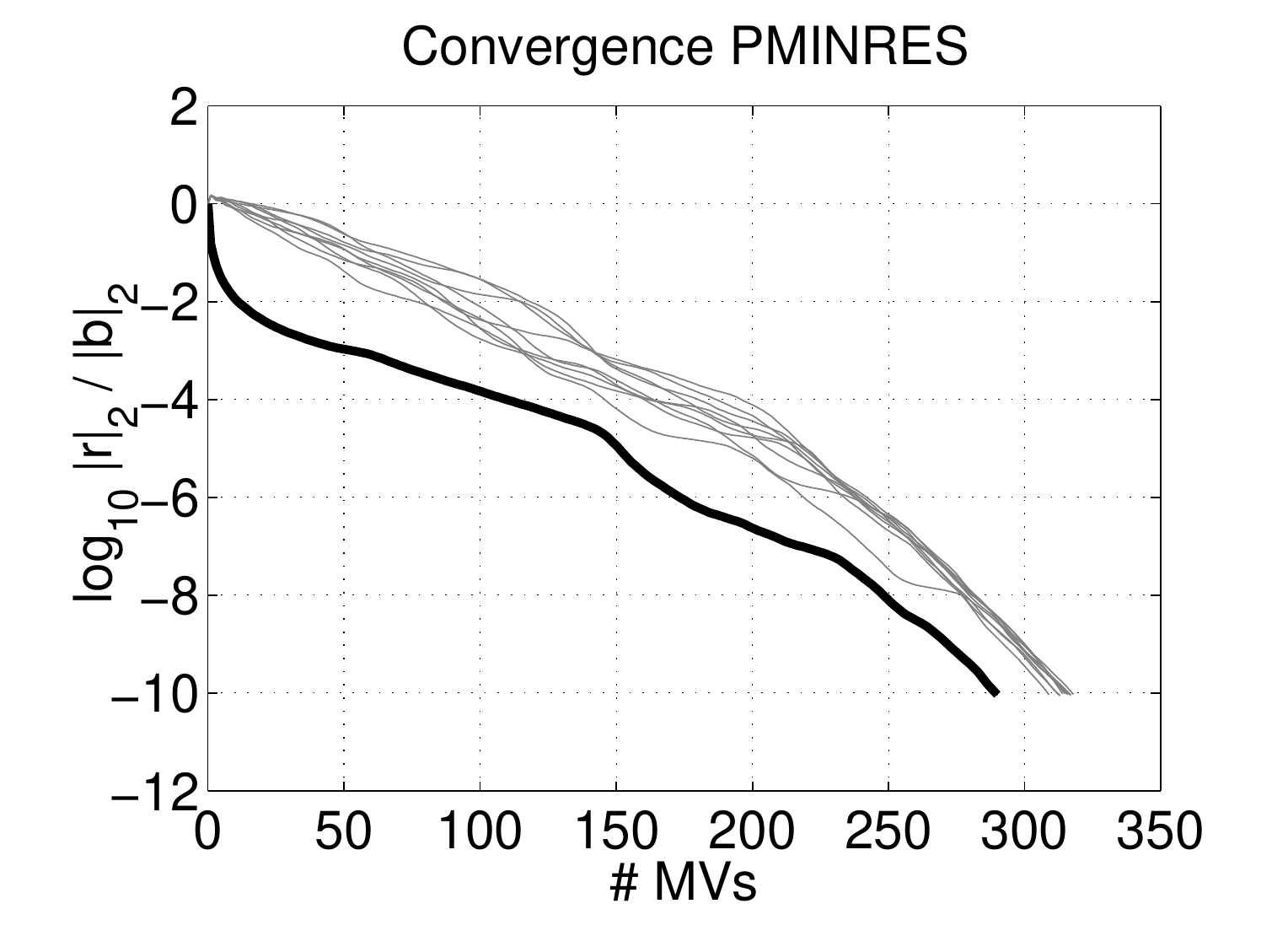}
	}
	\caption{Convergence for all \rhs-es of {sequence A} of POISSON for (a) \SRPCRap, (b) \PMINRES; thick black = first \rhs, thin gray = subsequent \rhs-es.}
	\label{fig:Poi_OldConv}
\end{figure}

\paragraph{Solving Sequence A}
We first solve for {sequence A} of $\bd$, cf. \eqref{eqn:RHSnatural}. For these \rhs-es the convergence curves of \SRPCRap and \MINRES are compared in figure~\ref{fig:Poi_OldConv}. We observe that \SRPCRap has a faster convergence after the recycling phase, possibly due to the space dimensions on which the residuals are orthogonalized. Anyway, for the first \rhs (from which the recycling space is built) the convergence speed-up by recycling is much different from that of the other \rhs-es.

\begin{figure}
	\subfigure[]{%
		\includegraphics[width=0.48\linewidth]{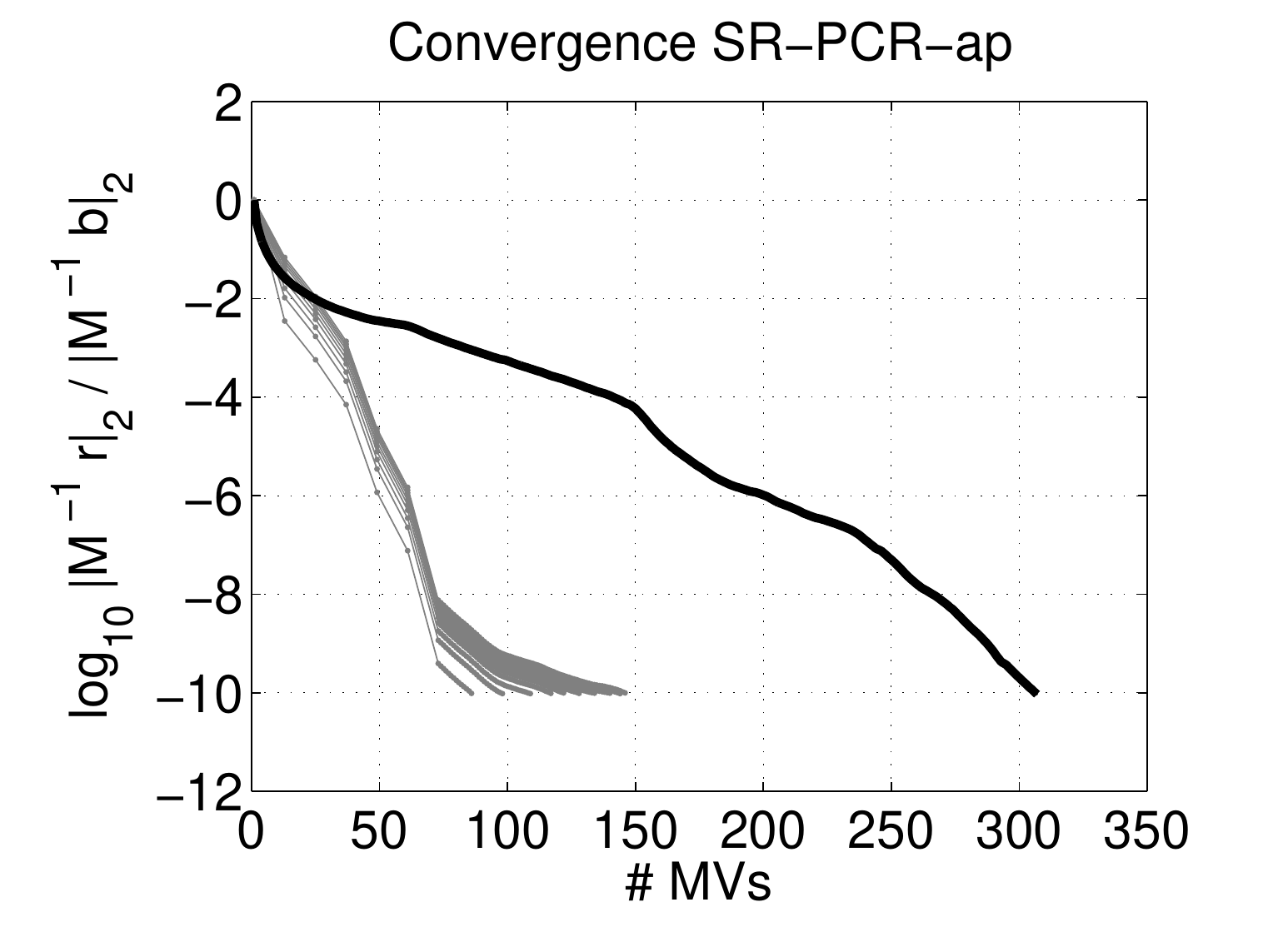}
	}
	\hfill
	\subfigure[]{%
		\includegraphics[width=0.48\linewidth]{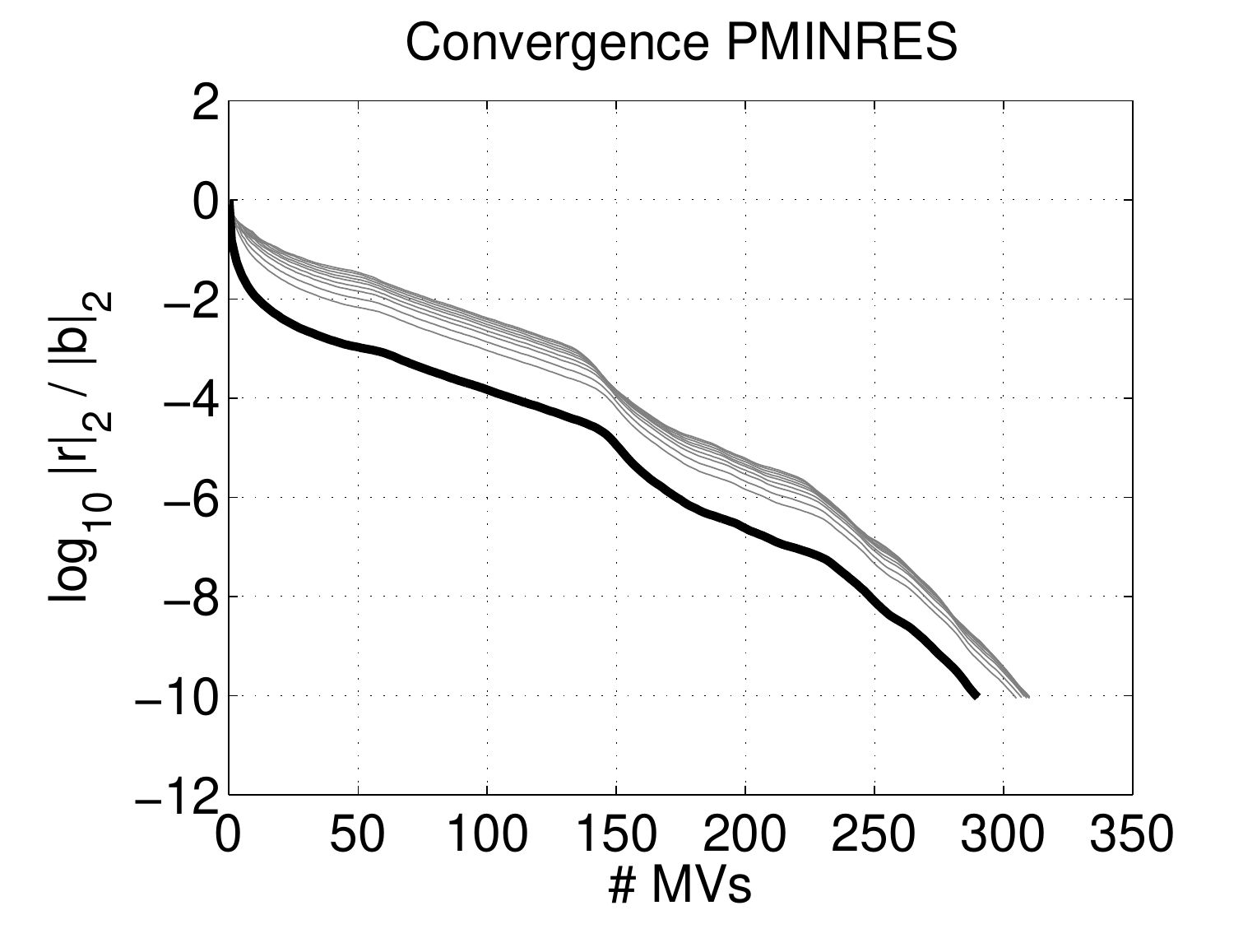}
	}
	\caption{Convergence for all \rhs-es of {sequence B} of POISSON for (a) \SRPCRap, (b) \PMINRES; thick black = first \rhs, thin gray = subsequent \rhs-es.}
	\label{fig:Poi_Conv}
\end{figure}

A possible reason for this might be that due to preconditioning the recycled search space and the \rhs-es do not fit together well. As already noted for sequence A, one can only expect that the the recycling space without preconditioning $\cK_m(\bA;\bd)$ (for $m$ sufficiently large) contains good candidates for subsequent \rhs-es. Instead, we used the recycling space $\bM^{-1}\,\cK_m(\bM^{-1}\,\bA;\bd)$.
\paragraph{Solving Sequence B}
As a second test we solve for {sequence B} of $\bd$, cf.  \eqref{eqn:RHSqnatural}. For this sequence we obtain completely different convergence curves than for {sequence A}, cf. figure~\ref{fig:Poi_Conv}. This is due to the fact that now the recycled search space contains useful candidates for subsequent \rhs-es for the preconditioned case. This can be directly seen in figure~\ref{fig:Poi_Conv} as the convergence curves do not differ very much from each other.

From the comparison of these two sequences we draw the following conclusion: \emph{If in a practical use case a sequence of type B with an arbitrary but known matrix $\bM$ occurs (such as e.g. a mass matrix from \FE discretizations), then for sake of a useful recycling space it might be advantageous to choose a preconditioner that is somehow \enquote{similar} to $\bM$.} This statement implies sequences of type A, as for these $\bM = \bI$ holds. So to a certain degree there is a trade-off between choosing an efficient versus a \enquote{similar} preconditioner when combining it with recycling.

\subsubsection{Stokes Problem}
We consider the following Stokes problem
\begin{align*}
\begin{array}{lllll}
- \Delta u 		& 				&{}+ \partial_x p &= 0 \quad &\text{in }\Omega\\
				&~- \Delta v 	&{}+ \partial_y p &= 0 \quad &\text{in }\Omega\\
\phantom{- }\partial_x u 	&{}+ \partial_y v & 				&= 0 \quad &\text{in }\Omega\\[1ex]
\phantom{- }n_x\nabla u &{}+ n_y \nabla v &{}- p \vec n &= 0 \quad &\text{on }\Gamma_N\\
&&\qquad u &= 0   \quad &\text{on }\Gamma_D\\
&&\qquad v &= g_D \quad &\text{on }\Gamma_D
\end{array}
\end{align*}
with $p, u, v$ denoting the pressure, velocity in $x$ and $y$ direction, respectively.
With the chosen boundary conditions this models flow in a lid-driven cavity.

\begin{figure}
	\subfigure[]{%
		\includegraphics[width=0.48\linewidth]{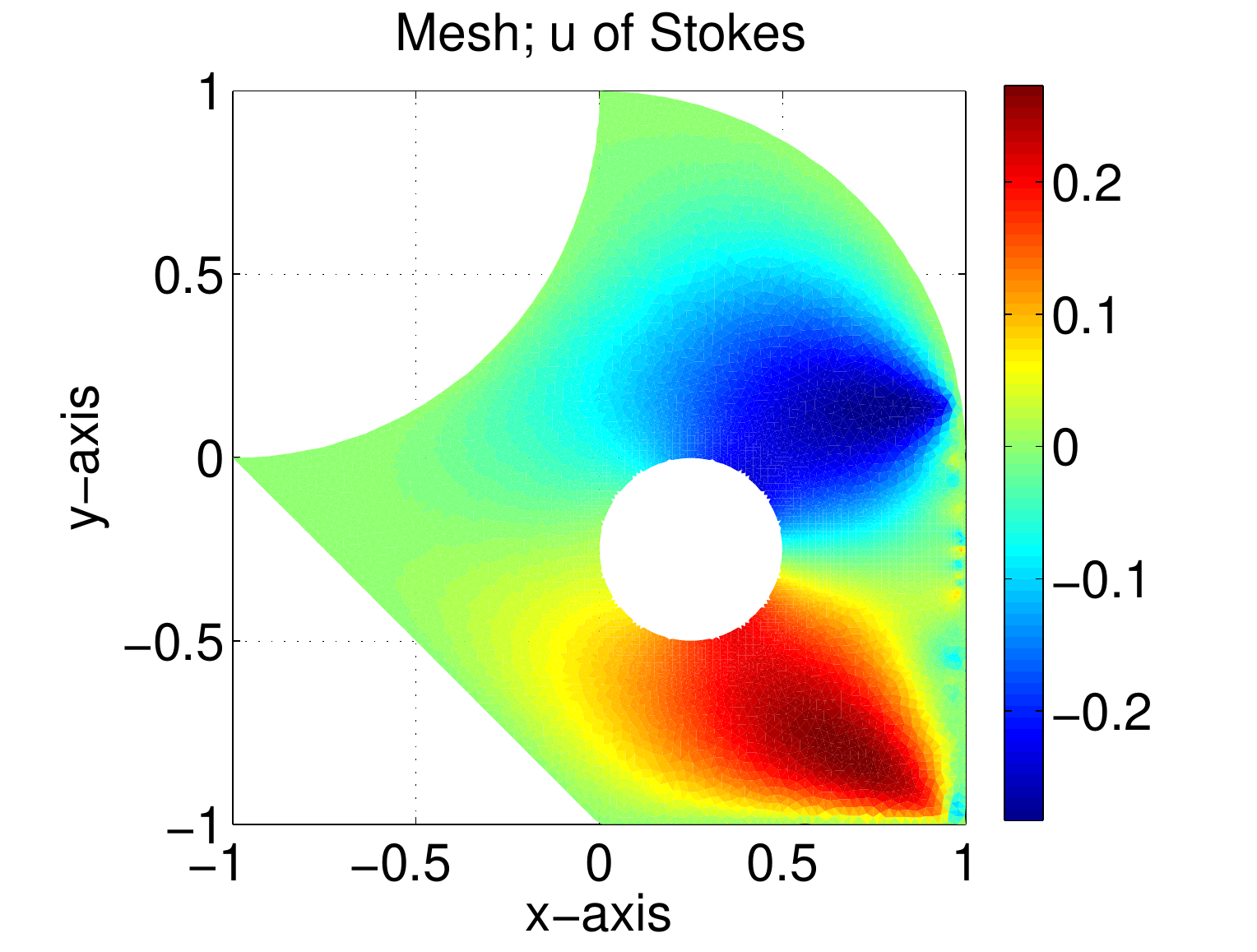}
	}
	\hfill
	\subfigure[]{%
		\includegraphics[width=0.48\linewidth]{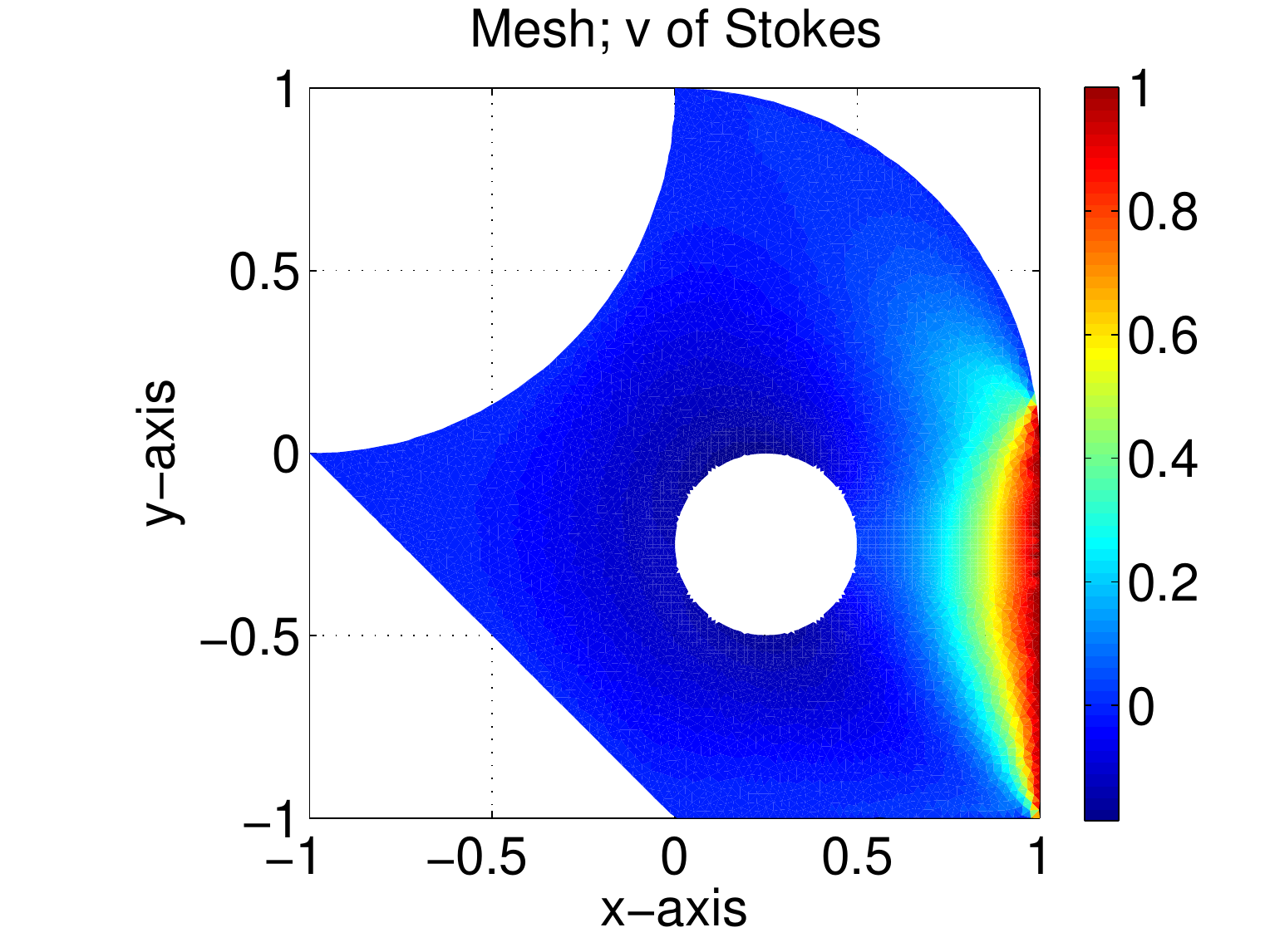}
	}
	\caption{Flow field of the numerical solution of the Stokes problem (a) in $x$- and (b) in $y$-direction.}
	\label{fig:Sto_Flo}
\end{figure}

\paragraph{Discretization}
We compute a numerical solution $u_h,v_h,p_h$ with $u_h,v_h$ piecewise quadratic \FE functions and $p_h$ piecewise linear \FE functions on the mesh in fig.~\ref{fig:Mesh1} (a).
By taking the discrete Laplacian $\bA$ from the Poisson problem above, $\bM_p \in \R^{6505 \times 6505}$ the \FE mass matrix for linear shape functions (index $_p$ indicates pressure space), and $\bB_i = \langle \psi, \partial_i \varphi \rangle_{L^2(\Omega)}$, $i \in \lbrace x,y \rbrace$, for quadratic shape functions $\varphi$ as columns and linear shape functions $\psi$ as rows, we solve
\begin{align*}
	\bK \, \bx \equiv \begin{bmatrix}
	\bA & \bOo & \bB^H_x \\
	\bOo & \bA & \bB^H_y \\
	\bB_x & \bB_y & \bOo
	\end{bmatrix}
	\begin{pmatrix}
	\bu \\
	\bv \\
	-\bp
	\end{pmatrix} = \bg \,.
\end{align*}
For the symmetric indefinite saddle point matrix $\bK \in \R^{58637 \times 58637}$ with $\kappa_1(\bK) \approx 3.8 \cdot 10^6$, we choose the preconditioner $\bM = \bL  \bL^H$ with $\bL= \diag(\bL_A,\bL_A,\bL_{M_p})$ where $\bL_A$ and $\bL_{M_p}$ are incomplete Cholesky decompositions of $\bA$ and $\bM_p$, respectively, using Matlab's \texttt{ichol} with threshold $\epsilon=10^{-4}$. The structure of $\bK$ and $\bL$ is given in fig.~\ref{fig:Mesh1} (b).
The numerical solution for the \rhs $\bg = (\bO^T,\bd^T,\bO^T)^T$ with $\bd$ as in the Poisson problem is shown in fig. \ref{fig:Sto_Flo}.

\paragraph{Solving Sequences B and C}
We keep $\ell = 2$ blocks of each $k=10$ vectors and $J=4$ to recycle the search space $\bM^{-1} \cdot \cK_{80}(\bA  \bM^{-1};\bb\ha)$. We emphasize that for the Stokes problem we use a strong preconditioner. To investigate its influence on the utility of the recycling space, we compare the results for sequence B and C with starting vector $\bg$, respectively. Note that by construction of both sequences the computed flow field is not discretely divergence-free as the third sub-vector of the \rhs-es is inconsistent.

\begin{figure}
	\subfigure[]{%
		\includegraphics[width=0.48\linewidth]{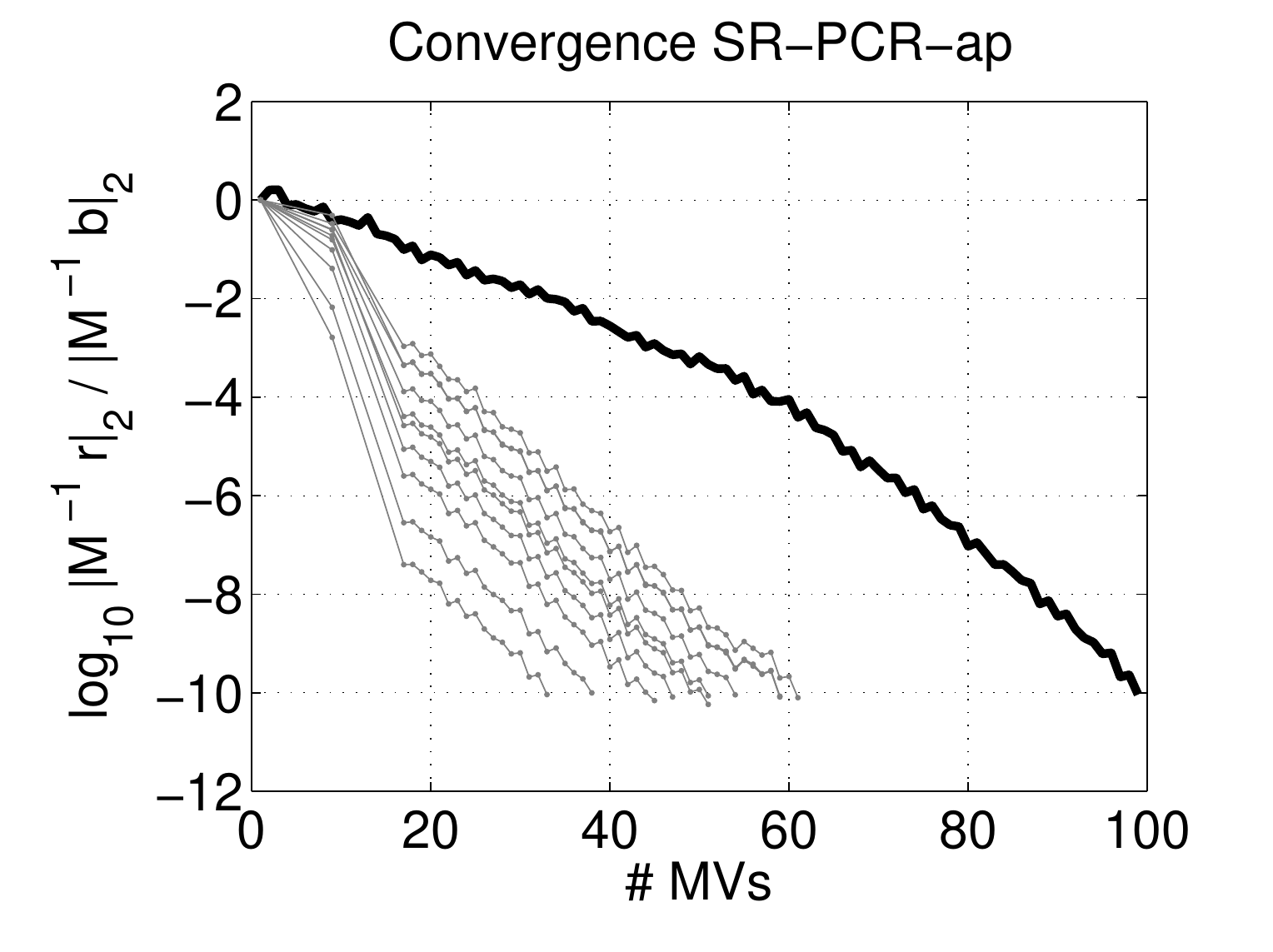}
	}
	\hfill
	\subfigure[]{%
		\includegraphics[width=0.48\linewidth]{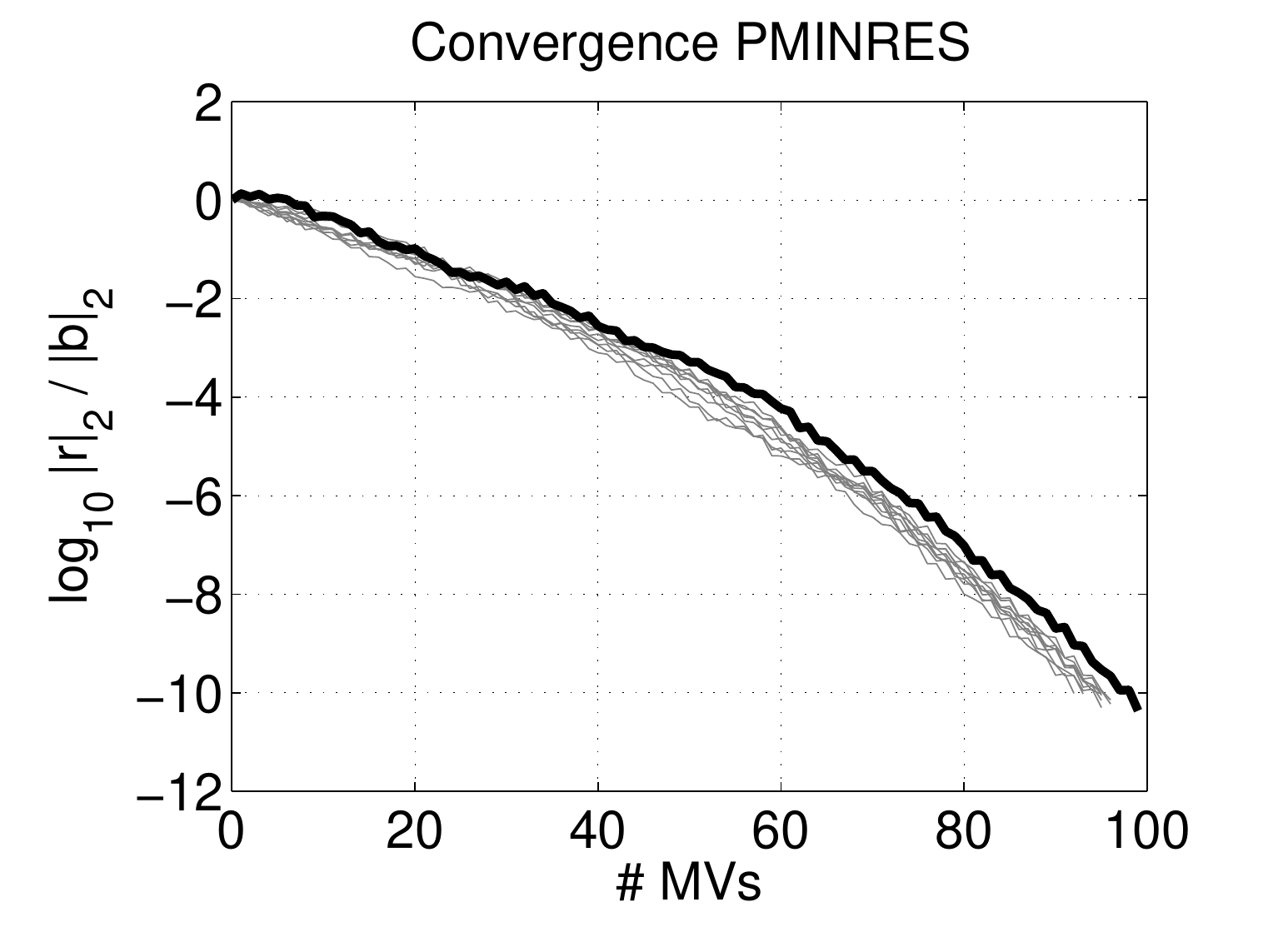}
	}
	\caption{Convergence for all \rhs-es of {sequence B} of STOKES for (a) \SRPCRap, (b) \PMINRES; thick black = first \rhs, thin gray = subsequent \rhs-es.}
	\label{fig:Sto_ConvB}
\end{figure}

\begin{figure}
	\subfigure[]{%
		\includegraphics[width=0.48\linewidth]{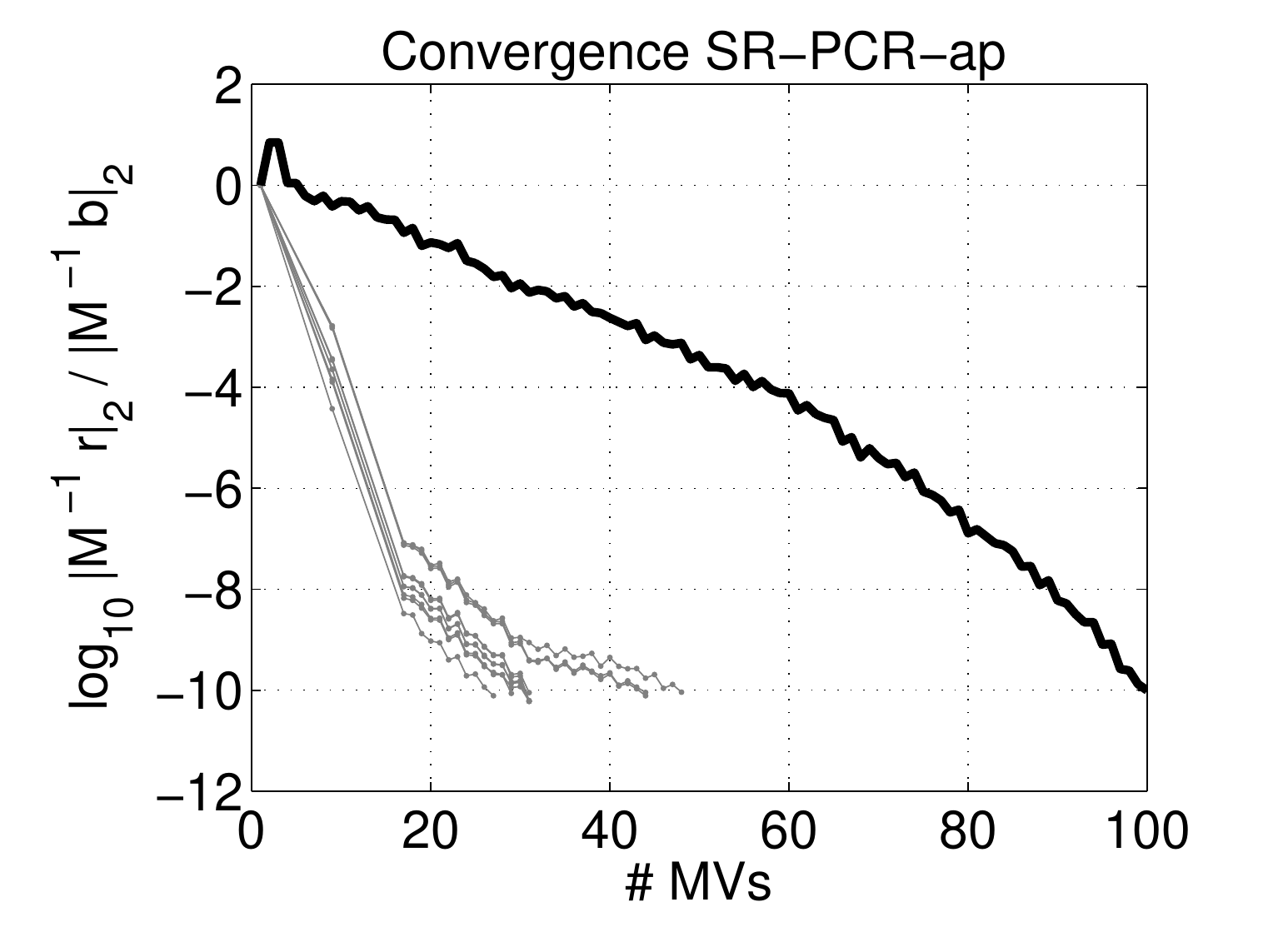}
	}
	\hfill
	\subfigure[]{%
		\includegraphics[width=0.48\linewidth]{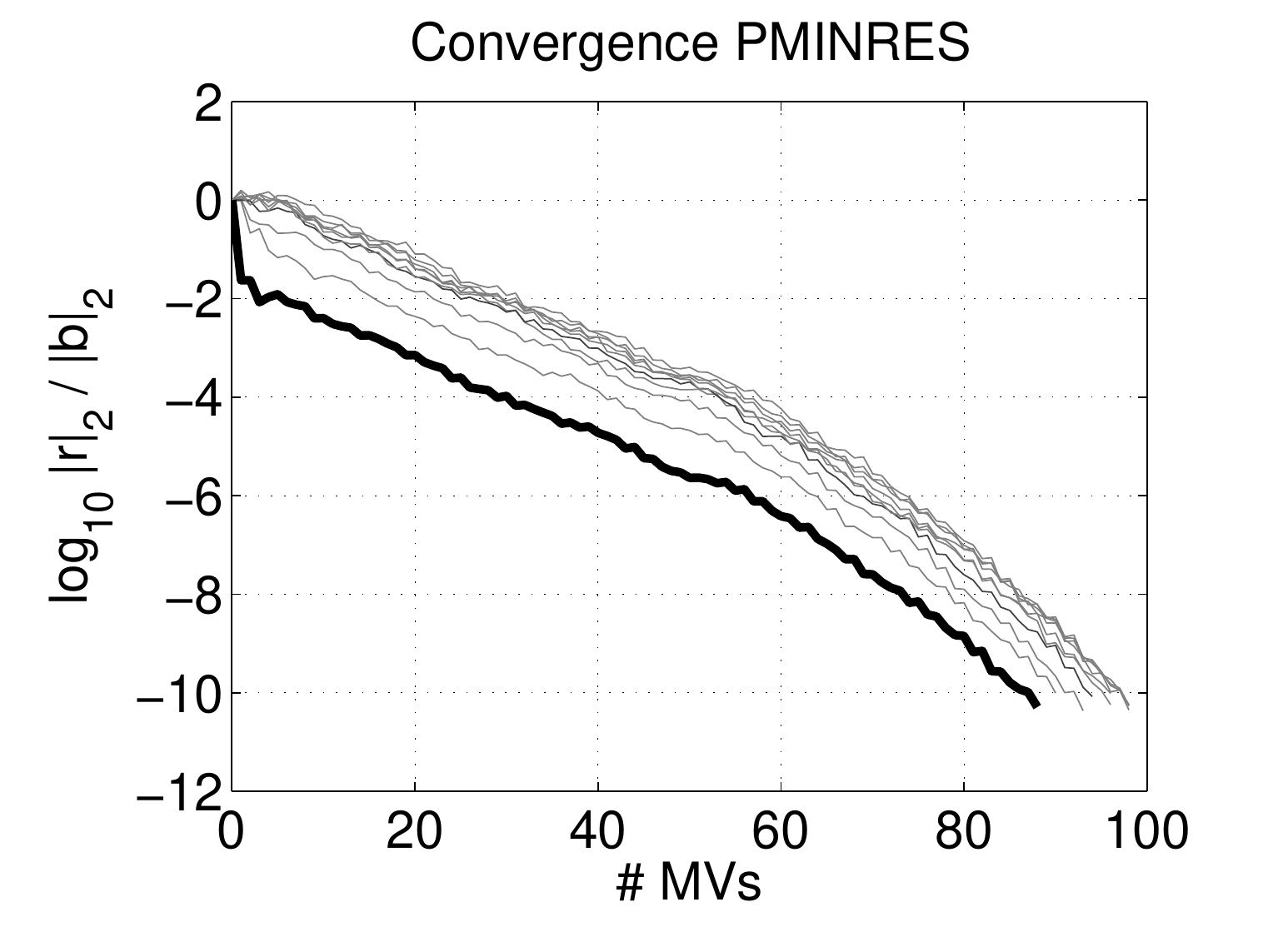}
	}
	\caption{Convergence for all \rhs-es of {sequence C} of STOKES for (a) \SRPCRap, (b) \PMINRES; thick black = first \rhs, thin gray = subsequent \rhs-es.}
	\label{fig:Sto_ConvC}
\end{figure}

The convergence curves for each sequence are given in fig.~\ref{fig:Sto_ConvB} and fig.~\ref{fig:Sto_ConvC}, respectively. From the figures we see that \SRPCRap offers a reduction in the average number of \MVec-s of about 50\% for both sequences compared to \PMINRES. However, for sequence B the convergence curves of \SRPCRap differ more from each other, with growing iteration numbers for later \rhs-es. This indicates that the recycling space gradually becomes out of date and that this effect is more pronounced for sequence B than for sequence C.

\FloatBarrier

\section{Conclusions and Outlook} \label{sec:conclusion}

The numerical experiments indicate that despite the delicate numerical properties of power and Horner schemes our method is quite practicable and efficient. In the experiments the number of stored columns could be reduced to $1/4^\text{th}$ of the dimension of the recycling space. The recycling approach shown here for \PCR can be easily adapted to other Krylov subspace methods.
From the experiments (cf. figs.~\ref{fig:Comp_CURL} and \ref{fig:Comp_CURLprec}) we deduce  that the proneness to rounding errors of \SRPCRap is correlated to the conditioning of the respective system.

Depending on the type of sequence, the recycling strategy turns out to be more or less successful, especially when combined with a preconditioner, cf. figs.~\ref{fig:Poi_OldConv} and \ref{fig:Poi_Conv}. As the preconditioner plays a dominant role in the construction of the recycling space it crucially influences its utility for subsequent systems. So in the design of a preconditioner the type of the \rhs-sequence might be taken into consideration when combining preconditioning with recycling. This is a current topic of on-going research.

\FloatBarrier
%\newpage


\begin{thebibliography}{99}

\footnotesize



%\bibitem[GL96]{MatComp} G. H. Golub and C. F. van Loan,
%	{Matrix Computations}, 3rd edition, John Hopkins University Press, 1996.


% Krylov Recycling Methods based on Compressing

%\bibitem[WA02]{Wang} S. Wang,
%	{Krylov subspace methods for topology optimization on adaptive meshes}, Ph.D. thesis, Tsinghua University, 2002.
	% R-MINRES

\bibitem{topol} S. Wang and E. de Sturler and G. H. Paulino,
	{Large-scale topology optimization using preconditioned Krylov subspace methods with recycling}, Int. J. for Num. Meth. in Engineering, Vol. 69(12), pp. 2441-2468, 2006.
	% R-MINRES

\bibitem{RBiCG} K. Ahuja and E. de Sturler and S. Gugercin and E. R. Chang,
	{Recycling BiCG with an Application to Model Reduction}, SIAM J. Sci. Comput. Vol. 34, No. 4, pp. A1925-A1949, 2012.
	% R-BiCG

\bibitem{MaxPlanck1} K. Ahuja and E. de Sturler and P. Benner,
	{Recycling BiCGSTAB with an Application to Parametric Model Order Reduction}, MPI Magdeburg preprints, pp. 13-21, 2013.
	% R-BiCGstab

\bibitem{BeFe} P. Benner and L. Feng,
	{Recycling Krylov Subspaces for Solving Linear	Systems with successively changing Right-Hand-Sides arising in Model Reduction}, Lecture Notes in Electrical Engineering, Vol. 74, pp. 125-140, Springer 2011.
	% R-GCR

%\bibitem[AH11]{AhujaDiss} K. Ahuja,
%	{Recycling Krylov Subspaces and Preconditioners}, Ph.D. thesis, Virginia Polytechnic Institute and State University, 2011.
	% R-BiCG

\bibitem{PSMJM06} M. Parks and E. de Sturler and G. Mackey and D.D. Johnson and S. Maiti,
	{Recycling Krylov subspaces for sequences of linear systems}, SIAM J. Sci. Comput. Vol. 28(5), pp. 1651-1674, 2006.
	% GCRO-DR

\bibitem{BBF13} J. Bolten and N. Bozovic and A. Frommer,
	{Preconditioning of Krylov subspace methods using recycling in Lattice QCD computations}, Proc. Appl. Math. Mech., Vol. 13, pp. 413-414, 2013.


% Applications
\bibitem{NavStokes} K. Mohamed and S. Nadarajah and M. Paraschivoiu,
	{Krylov Recycling techniques for unsteady simulation of turbulent aerodynamic flows}, 26th international congress of the aeronautical sciences, 2008.

\bibitem{Circuit} Z. Ye and Z. Zhu and J. R. Phillips,
	{Generalized Krylov Recycling Methods for Solution of Multiple Related Linear Equation Systems in Electromagnetic Analysis}, Design Automation Conference 2008, p. 682-687.

%\bibitem[KS06]{KiSt} M. Kilmer and E. de Sturler,
%	{Recycling Subspace Information for diffuse optical Tomography}, SIAM J. Sci. Comput., Vol. 27(6), pp. 2140-2166, 2006.

\bibitem{SSX13} K. M. Soodhalter and D. E. Szyld and F. Xue,
	{Krylov Subspace Recycling for Sequences of Shifted Linear Systems}, Elsvier J. Appl. Num. Math. Vol. 81, pp. 105-118, 2014.


% Restarting
\bibitem{Morgan2} R. B. Morgan,
	{A restarted GMRES method augmented with eigenvectors}. SIAM J. Matrix Anal. Appl., 16:1154-1171, 1995.

\bibitem{GMRESR} H. A. van der Vorst and C. Vuik,
	{GMRESR: A family of nested GMRES methods.} Num. Lin. Alg. with Appl., 1:369-386, 1994.

\bibitem{Morgan3} R. B. Morgan,
	{Implicitly restarted GMRES and Arnoldi methods for nonsymmetric systems of equations}. SIAM J. Matrix Anal. Appl., 21 (2000), pp. 1112-1135.


\bibitem{HaNo} M. Habu and T. Nodera,
	{GMRES$M$ algorithm with changing the restart cycle adaptively}, Proceedings of ALGORITHMY 2000, Conference on Scientific Computing, pp. 254-263.





% Compressing Approaches
\bibitem{GCROT} E. de Sturler,
	{Truncation Strategies for optimal Krylov subspace methods}, SIAM J. Numer. Anal., Vol. 36(3), pp. 864-889, 1999.
	% Optimal Truncation

\bibitem{Morgan1} R. B. Morgan,
	{GMRES with Deflated Restarting}, SIAM J. Sci. Comput., 24(1), pp. 20-37, 2002.
	% Deflation


% Short Representations
\bibitem{MyReport} M. P. Neuenhofen,
	{Short-Recurrence and -Storage Recycling of large Krylov-Subspaces for Sequences of Linear Systems with changing Right-Hand-Sides}, Technical Report, available on arXiv: 1512.05101, 2015.


% \bibitem[NE15]{MySlides} M. P. Neuenhofen,
% 	{Short Recycling of Krylov Subspaces}. Talk for NA group, TU Delft. \texttt{MartinNeuenhofen.de/ShortRecycling/Short\_Recycling.html}, 2015.
%

\bibitem{SimpleFEM} M. P. Neuenhofen,
	{Simple Finite Element Code},\\ \texttt{MartinNeuenhofen.de/Simple\_FEM/SimpleFEM.html}, 2015.


% Rounding on Lanczos
\bibitem{LancPartReortho} H. D. Simon,
	{The Lanczos Algorithm With Partial Reorthogonalization}, Mathematics of Computation, 46(165):115-142, 1984.


\bibitem{PaiSau75} C. C. Paige and M. A. Saunders,
	{Solution of sparse indefinite Systems of linear equations}, SIAM J. Numer. Anal. Vol. 12(4), pp. 617-629, 1975.


% Singular Systems
\bibitem{MRQLP} S. T. Choi, C. C. Paige and M. A. Saunders, {MINRES-QLP: A Krylov Subspace Method for Indefinite or Singular Symmetric Systems}, SIAM J. Sci. Comput., 33(4), pp. 1810-1836, 2011.

\bibitem{PCRsing} K. Hayami, {On the Behaviour of the Conjugate Residual Method for Singular Systems}, Technical Report, National Institute of Informatics, Tokyo 2011.


% Testcases
\bibitem{FSMC} T. A. Davis and Y. Hu,
	{The University of Florida Sparse Matrix Collection}. ACM Transactions on Mathematical Software, 38(1):1-25. \texttt{http://www.cise.ufl.edu/research/sparse/matrices}. 2011.


\bibitem{IDRstab} G. L. G. Sleijpen and M. B. van Gijzen,
	{Exploiting BiCGstab($\ell$) Strategies to Induce Dimension Reduction}, SIAM J. Sci. Comput. 32(5):2687-2709, 2010.


\bibitem{Boeing} Duff, I. S. and R. G. Grimes and J. G. Lewis,
	{Sparse Matrix Problems}. ACM Trans. on Mathematical Software, 14(1):1-14, 1989.


\bibitem{GR11} S. Gross and A. Reusken,
	{Numerical Methods for Two-phase Incompressible Flows}. First edition 2011, Springer Series in Computational Mathematics, Vol. 40, ISBN 978-3-642-19685-0.

\bibitem{CURL} Z. Bai,
	{Krylov subspace techniques for reduced-order modeling of large-scale dynamical systems}. {Applied Numerical Mathematics},	43(1-2):9-44, 2002.

% Books
\bibitem{Saad1} Y. Saad, {Iterative Methods for Sparse Linear Systems, 2nd edition}, SIAM, 2000.

\small


\end{thebibliography}
\end{document}